\documentclass[10pt, twoside]{amsart}
\usepackage{amsmath,amsthm,amssymb}
\usepackage{times}
\usepackage{enumerate}

\usepackage[letterpaper, top=2.9cm, bottom=2.9cm,left=2.8cm, right=2.8cm, heightrounded,bindingoffset=0mm]{geometry}

\theoremstyle{definition}
\newtheorem{defi}{Definition}[subsection]

\newtheorem{rmk}[defi]{Remark}
\newtheorem{eg}[defi]{Example}

\theoremstyle{plain}
\newtheorem{thm}[defi]{Theorem}
\newtheorem{lem}[defi]{Lemma}
\newtheorem{prop}[defi]{Proposition}
\newtheorem{cor}[defi]{Corollary}

\newtheorem*{claim}{Claim}
\newtheorem*{theoremint}{Theorem}
\newtheorem*{propint}{Proposition}
\newtheorem*{corint}{Corollary}

\theoremstyle{remark}
\newtheorem*{notation}{Notation}

\newtheorem*{setting}{Setting and notation}

\numberwithin{equation}{section}

\usepackage{fancyhdr}

\usepackage{amsfonts}
\usepackage{amscd}
\usepackage{mathrsfs}

\usepackage[all,cmtip]{xy}	
\xyoption{arrow}

\usepackage{color}	

\usepackage{graphicx}	

\usepackage[english]{babel} 
\usepackage[utf8]{inputenc}
\usepackage[T1]{fontenc}

\usepackage{enumitem}

\newcommand{\NN}{\mathbb{N}}
\newcommand{\ZZ}{\mathbb{Z}}

\newcommand{\PP}{\mathbb{P}}

\newcommand{\VV}{\mathbb{V}}

\newcommand{\A}{\mathcal{A}}
\newcommand{\B}{\mathcal{B}}
\newcommand{\V}{\mathcal{V}}
\newcommand{\W}{\mathcal{W}}
\newcommand{\Z}{\mathcal{Z}}

\newcommand{\Ff}{{\mathcal{F}}}

\newcommand{\Xx}{\mathcal{X}}

\newcommand{\E}{\mathscr{E}}

\newcommand{\Oo}{\mathscr{O}}

\newcommand{\p}{{\mathfrak{p}}}

\newcommand{\Rr}{\mathcal{R}}
\newcommand{\Rrt}{{\widetilde{\mathcal{R}}}}
\newcommand{\Xt}{{\widetilde{X}}}

\newcommand{\Xxt}{{\widetilde{\Xx}}}
\newcommand{\xt}{{\widetilde{x}}}
\newcommand{\Cun}{{X_{univ}}}
\newcommand{\Ctun}{{\widetilde{X}_{univ}}}

\newcommand{\Op}{\widehat{\mathcal{O}}}
\newcommand{\Ll}{{\mathcal{L}}}
\newcommand{\lL}{{\mathfrak{l}}}
\newcommand{\Aa}{{\mathcal{A}}}

\newcommand{\Hur}[1][1]{{\mathcal{H}\textrm{ur}(\Gamma,\xi)_{g, #1}}}
\newcommand{\Hurz}{{\mathcal{H}\textrm{ur}(\Gamma,\xi)_{g}}}
\newcommand{\bHurz}{\overline{{\mathcal{H}\textrm{ur}}}(\Gamma,\xi)_{g}}
\newcommand{\bHur}[1][1]{{\overline{\mathcal{H}\textrm{ur}}(\Gamma,\xi)_{g, #1}}}
\newcommand{\bHurs}[2]{{\overline{\mathcal{H}\textrm{ur}}(\Gamma,#2)^{#1}_{g}}}
\newcommand{\Hurs}[2]{{\mathcal{H}\textrm{ur}(\Gamma,#2)^{#1}_{g}}}

\newcommand{\Mg}{{\mathcal{M}_{g}}}

\newcommand{\Mgn}[1][1]{{\mathcal{M}_{g, #1}}}
\newcommand{\bMgn}[1][1]{{\overline{\mathcal{M}}_{g, #1}}}

\newcommand{\Hom}{\textrm{Hom}}
\newcommand{\Ext}{\textrm{Ext}}

\newcommand{\Extt}{\mathscr{E}\!\textit{xt}}

\newcommand{\Aut}{\text{Aut}}
\newcommand{\End}{\textrm{End}}
\newcommand{\Def}{\textrm{Def}}

\newcommand{\Res}{\textrm{Res}}
\newcommand{\Sym}{\textrm{Sym}}
\newcommand{\Char}{\textrm{Char}}
\newcommand{\Spec}{\textrm{Spec}}
\newcommand{\SPec}{\textbf{Spec}}

\newcommand{\Bun}{\textrm{Bun}}

\newcommand{\Lie}{\text{Lie}}

\newcommand{\Repl}{{\text{IrRep}_\ell}}
\newcommand{\AND}{\quad \text{and} \quad}
\newcommand{\univ}{{\text{univ}}}

\newcommand{\Qdx}{{\, \big/}}
\newcommand{\Qsx}{{\, \big\backslash}}
\newcommand{\circdot}{{_\circ^\circ}}

\newcommand{\g}{{\mathfrak{g}}}
\newcommand{\h}{{\mathfrak{h}}}

\newcommand{\hA}{{\mathfrak{h}_{\mathcal{A}}}}
\newcommand{\hB}{{\mathfrak{h}_{\mathcal{B}}}}
\newcommand{\hL}{{\mathfrak{h}_{\mathcal{L}}}}
\newcommand{\hLi}[1][i]{{\mathfrak{h}_{\mathcal{L}_{#1}}}}

\newcommand{\HhL}{{\widehat{\mathfrak{h}}_{\mathcal{L}}}}
\newcommand{\HhLi}[1][i]{{\widehat{\mathfrak{h}}_{\mathcal{L}_{#1}}}}

\newcommand{\gL}{{\mathfrak{g}{\mathcal{L}}}}
\newcommand{\HgL}{{\widehat{\mathfrak{g}}{\mathcal{L}}}}

\newcommand{\Hhn}[1]{{\mathcal{H}_\ell({#1})}}
\newcommand{\tHh}[1]{{\widetilde{\mathcal{H}_\ell}({#1})}}
\newcommand{\Hhne}{{\mathcal{H}_\ell(0)}}
\newcommand{\tHhe}{{\widetilde{\mathcal{H}_\ell}({0})}}

\newcommand{\tm}{t_-}
\newcommand{\tp}{t_+}
\newcommand{\tmj}[1][j]{t_{#1,-}}
\newcommand{\tpj}[1][j]{t_{#1,+}}
\newcommand{\tpm}{t_\pm}
\newcommand{\tpmj}[1][j]{t_{#1,\pm}}

\newcommand{\SLn}[1][n]{\text{SL}_{#1}}

\newcommand{\gl}{{\mathfrak{gl}}}

\newcommand{\T}{\mathcal{T}}
\renewcommand{\H}{\mathcal{H}}

\usepackage[colorlinks=true]{hyperref}
\hypersetup{bookmarksnumbered=true, linkcolor=black, citecolor=blue, urlcolor=black,}

\begin{document}

\normalsize

\title{Conformal blocks attached to twisted groups}

\author{Chiara Damiolini}

\email{chiara.damiolinid@gmail.com}

\maketitle

\begin{abstract} The aim of this paper is to generalize the notion of conformal blocks to the situation in which the Lie algebra they are attached to is replaced with a sheaf of Lie algebras depending on covering data of curves. The result is a vector bundle of finite rank on the stack $\bHur[n]$ parametrizing $\Gamma$-coverings of curves. Many features of the classical sheaves of conformal blocks are proved to hold in this more general setting, in particular the factorization rules, the propagation of vacua and the WZW connection.
\end{abstract}

\section{Introduction}
In conformal field theory \cite{tsuchiya1989conformal} there is a way to associate to a natural number $\ell$ and to a simple and simply connected group $G$ over an algebraically closed field $k$ of characteristic zero, a vector bundle $\VV_\ell(0)$, called the \emph{sheaf of covacua}, on $\Mg$, the stack parametrizing smooth curves of genus $g$. The goal of this paper is to generalize the construction and properties of these bundles to the case in which the group $G$ they are attached to is replaces by a \emph{parahoric Bruhat-Tits group} $\H$ which depends on cyclic coverings of curves. For this reason, this new sheaf of covacua will not be defined on $\Mg$, but on the Hurwitz stack parametrizing coverings of curves. 

One of the reasons that lead algebraic geometers to study the sheaves of covacua, and generalization thereof, is their relation to the stack parametrizing principal bundles on curves. In the classical setting, the fibres of the dual of $\VV_\ell(0)$ over the curve $X$, called space of \emph{conformal blocks}, have been identified as global sections of well determined line bundles on $\Bun_{G,X}$, obtaining insights on the geometry of this moduli space  \cite{BeauvilleLaszlo1994Conformal} \cite{Faltings1994Verlinde}\cite{LS1997PicardBunG} \cite{sorger1996formule}. The key point to prove this isomorphism is the \emph{uniformization theorem} which describes $\Bun_{G,X}$ as a quotient of the the affine Grassmannian $\textsf{Gr}(G)$, whose Picard group and the space of global sections of line bundles have been described in terms of representations of $\g$ by Kumar \cite{Kumar} and Mathieu \cite{Mathieu}. The uniformization theorem, which was proved initially by Beauville and Laszlo in \cite{BeauvilleLaszlo1994Conformal} for $G=\SLn$, has been shown to hold for any simple group $G$ by Drinfeld and Simpson \cite{DrinfeldSimpson}. The identification of global sections of line bundles and conformal blocks was then extended to the case of parabolic bundles by Pauly in \cite{Pauly1996Parabolic} and by Laszlo and Sorger in \cite{LS1997PicardBunG}. Finally Heinloth, answering questions posed by Pappas and Rapoport in \cite{PappasRapoport2007Questions}, proved the uniformization theorem for $\Bun_{\H, X}$, for connected parahoric Bruhat-Tits groups $\H$ in \cite{heinloth2010uniformization}, where he also gave a description of the Picard group of $\Bun_{\H,X}$. 

This motivates out interest in developing the theory of conformal blocks attached to twisted groups $\H$. Inspired by \cite{balaji2011moduli} we restrict ourselves to consider only those groups \emph{arising from cyclic Galois coverings}, so that $\Mg$ will be replaced by the Hurwitz stack $\Hurz$. 

As in the classical case, also in this new setting the sheaf of conformal blocks arises as dual of the \emph{sheaf of covacua}. Classically, for every $\ell \in \ZZ_{> 0}$, given $n$ dominant weights $\lambda_i$ of level at most $\ell$ of $\g=\Lie(G)$, it is possible to construct vector bundles $\VV_\ell(\lambda_1, \dots, \lambda_n)$ on $\bMgn[n]$, called sheaves of covacua \cite{tsuchiya1989conformal}. Similarly, in Section \ref{sec-ConformalBlocks} we show how to associate to admissible representations $\V_1, \dots, \V_n$ of fibres of $\h=\Lie(\H)$ outside the branch locus of the covering, a sheaf $\VV_\ell(\V_1, \dots, \V_n)$ on $\bHur[n]$ which generalizes the definition and the properties of the classical covacua. 

Besides the relation with $\Bun_{G,X}$, the classical sheaves of conformal blocks have been used to investigate the geometry of $\bMgn[n]$ and provide examples of divisors satisfying interesting combinatorial relations. In particular, in the case $g=0$, the sheaves of covacua are globally generated, hence their first Chern classes, called \emph{conformal blocks divisors}, lie in the nef cone, where they generate a full dimensional sub-cone \cite{Fak12CBdivisors}. The combinatorial nature of these divisors, inherited from the structure of the sheaf of covacua, made possible to establish vanishing and non-vanishing criteria for conformal blocks divisors which gave insights in the geometry of $\bMgn[n]$ \cite{BGM15Vanishing} \cite{BGM16Nonvanishing}. In similar way, the twisted sheaves of covacua and conformal blocks provide vector bundles with a rich combinatorial structure on $\bHur[n]$ which can be used to study $\bHur[n]$ itself and, and forgetting about the covering data, the geometry of moduli of curves.\\


The main results of this paper (Theorem \ref{thm-LocFree} and Corollary \ref{cor-ConnectionLocFree}) can be resumed in the following statement.

\begin{theoremint} The sheaf $\VV_\ell(\V_1, \dots, \V_n)$ is a vector bundle of finite rank on $\bHur[n]$ which admits a projectively flat connection on $\Hur[n]$.
\end{theoremint}

In Section \ref{sec-factorizationRulesS} we describe the properties of these sheaves. In particular Proposition \ref{prop-PropVacua} generalizes the so called \emph{propagation of vacua}, which essentially says that trivial representations $\V(0)$ do not modify the sheaf of conformal block. As a consequence we have the following result (Corollary \ref{cor-PropVacua}.)

\begin{propint} The bundle $\VV_\ell(\V(0), \dots, \V(0))$ is independent of the choice of the marked points, hence it descends to a vector bundle $\VV_\ell(0)$ on $\Hurz$.
\end{propint}

As in the classical case, in Proposition \ref{prop-NodalDegenerationH} we formulate \emph{factorization rules} controlling the rank of the vector bundle under degeneration of the covering. This is the key input to use induction on the genus of the curves to achieve a formula computing the rank of the bundles.

\begin{propint} Let $(q \colon \Xt \to X, \p) \in \bHur(k)$ such that $X$ is irreducible and has only one nodal point $x$. Let $X_N$ be its normalization so that $q_N \colon \Xt_N \to X_N$ is a $\Gamma$-covering with three marked points. Then for any $\W \in \Repl(\h|_\p)$ we have a canonical isomorphism
\[ \VV_\ell(\W)_X \cong \bigoplus_{\V \in \Repl(\h|_x)} \VV_\ell(\W, \V, \V^*)_{X_N}.
\]
\end{propint}

In the untwisted case this factorization property, shown in \cite{tsuchiya1989conformal}, allowed to reduce the computation of the rank of the sheaves of conformal blocks to the case of $\PP^1$ with three marked points and obtaining in this was the \emph{Verlinde formula} to compute dimension of line bundles on $\Bun_{G,X}$  \cite{Faltings1994Verlinde} \cite{sorger1996formule}.\\

It is important to remark that after the first draft of this paper was completed, several authors worked towards a definition and properties of conformal blocks attached to twisted groups arising from coverings. In particular Zelaci in \cite{zelaci2017moduli} constructs special cases of twisted conformal blocks and relates them to global sections of line bundles on appropriate $\Bun_\H$. In the more recent pre-print \cite{KumarHong}, Hong and Kumar construct conformal blocks attached to groups arising from Galois coverings which are not necessarily cyclic, generalizing the factorization rules and projective connection, and establishing their relation to global sections of line bundles. It is worth to mention that their construction is compatible with the one presented in \cite{zelaci2017moduli} and in this paper.\\

We now give an overview of how the twisted conformal blocks are defined, generalizing the methods used in \cite{kac1994infinite}, \cite{tsuchiya1989conformal} and \cite{looijenga2005conformal}. To start with, let us briefly explain how a Galois covering determines a twisted group. We fix the cyclic group $\Gamma:=\ZZ/p\ZZ$ of prime order $p$ and a group homomorphism $\rho \colon \Gamma \to \Aut(G)$. To every $\Gamma$-covering $q \colon \Xt \to X$ of curves, we associate the sheaf of groups $\H$ defined as the group of $\Gamma$-invariants of the Weil restriction of $\Xt \times_k G$ along $q$, i.e. $\H= q_*(\Xt \times_k G)^\Gamma$. The Lie algebra of $\H$ is denoted $\h$. The preliminaries on coverings, on the stack parametrizing them and on the construction of $\H$ are collected in Section \ref{sec-HurwitzStacks}.

To construct the sheaf of covacua, we explain how to associate to each Galois covering $(q \colon \Xt \to X, \p) \in \Hur(\Spec(k))$ and representation $\V$ of $\h|_{\p}$ a finite dimensional vector space. We want to remark that in our construction, we assume the point $\p$ marking $X$ is disjoint from the branch locus of $q$: this implies that $\h|_\p$ is isomorphic, although non canonically, to $\g$. It follows that once we choose such an isomorphism, we can use the classical construction \cite[Chapter 7]{kac1994infinite} to associate to each representation $\V \in \Repl(\h|_\p)$ the integrable highest weight representation $\Hhn{\V}$ of $\HhL$, a central extension of $\hL=\h|_{k(\!(t)\!)}$ defined in terms of Killing form and residue pairing. The key point is to see that this construction is actually independent of the isomorphism chosen between $\h|_\p$ and $\g$. Thanks to the residue theorem, the Lie algebra $\hA:= \h|_{X \setminus \p}$ is a Lie subalgebra of $\HhL$ and we set $\VV_\ell(\V)_X$ to be the quotient $\hA \Qsx \Hhn{\V}$. The construction of the sheaf of covacua runs similarly for any family of coverings $(\Xt \to X, \sigma) \in \bHur(S)$, being careful that the isomorphism between $\h|_{\sigma(S)}$ and $\g \otimes_k S$ exists only locally on $S$. The construction of the sheaf is the content of Section \ref{sec-ConformalBlocks}.

Although it is easy to show that $\VV_\ell(\V)$ is coherent (Proposition \ref{prop-HhnFinGen}), it is not immediate from its construction that it is also locally free. Following the approach of Looijenga in \cite{looijenga2005conformal}, the first step to achieve this result is to generalize to this twisted setting the \emph{Wess-Zumino-Witten connection} defined in terms of conformal field theory. After recalling how the connection arises using the Virasoro algebra of $\hL$, in Section \ref{sec-ProjConn} we show:

\begin{corint} [Corollary \ref{cor-ConnectionLocFree}] The sheaf $\VV_\ell(\V_1, \dots, \V_n)$ on $\bHur[n]$ is equipped with a projectively flat connection with logarithmic singularities along the boundary $\bHur[n] \setminus \Hur[n]$. \end{corint}

This shows in particular that $\VV_\ell(\V_1, \dots, \V_n)$ is a locally free module over $\Hur[n]$. Combining this with a refined version of the factorization rules (Proposition \ref{EisIso}), we are able to prove the local freeness of the sheaf on the whole stack $\bHur[n]$. Also in this twisted setting then, the factorization rules play a double role in the theory of conformal blocks. On one side they contribute to show that $\VV_\ell(\V_1, \dots, \V_n)$ is locally free on the whole $\bHur[n]$, and on the other side they are a useful tool to reduce the computation to lower genera curves.

\begin{setting} \label{subsec-NotationAndSetting}
Throughout the paper we fix the following objects.
\begin{itemize}
\item An algebraically closed field $k$ of characteristic zero.
\item A simple and simply connected algebraic group $G$ over $\Spec(k)$.
\item A prime $p$ and for simplicity of notation we denote the group $\ZZ/p\ZZ$ by $\Gamma$.
\item A group homomorphism $\rho \colon \Gamma \to \Aut(G)$.
\end{itemize}
\end{setting}


\section{Preliminaries on groups arising from coverings and Hurwitz stacks} \label{sec-HurwitzStacks}

In this section we introduce the group schemes associated with coverings as indicated in the introduction. Since we need to work with these groups in families, we will formulate the definition for families of coverings of curves. We obtain in this way the family $\H_{\univ}$ over the universal curve $\Cun$ over the Hurwitz stack parametrizing coverings of curves.

\begin{defi} Let $\pi \colon X \to S$ be a possibly nodal curve over a $k$-scheme $S$. A \emph{Galois covering} of $X$ with group $\Gamma$, called also \emph{$\Gamma$-covering}, is the data of 
\begin{enumerate}[label=\alph*)]
\item a finite, faithfully flat and generically étale map $q \colon \widetilde{X} \to X$ between curves;
\item an isomorphism $\phi \colon \Gamma \cong \Aut_X(\Xt)$;\end{enumerate}
satisfying the following conditions: \begin{enumerate}
\item each fibre of $\Xt$ is a generically étale $\Gamma$-torsor over $X$;
\item the singular locus of $\pi q$, i.e. the set of nodes of $\Xt$, is contained in the étale locus of $q$.
\end{enumerate}
\end{defi}

We want to attach to any $\Gamma$-covering $(\Xt \overset{q}\to X \overset{\pi} \to S)$ and to the homomorphism $\rho \colon \Gamma \to \Aut(G)$ a group scheme $\H$ over $X$ in the same fashion as in \cite[Section 4]{balaji2011moduli}. 

\begin{rmk}We remark that Balaji and Seshadri consider $\rho$ to map to the inner automorphisms of $G$ only, i.e. arising from a morphism $\Gamma \to G$. Without imposing that restriction we allow also groups $\H$ which are non-split over the generic point of $X$.\end{rmk} 

First of all we consider the scheme $\widetilde{G}:= \Xt \times_k G$ and let $q_*(\widetilde{G})$ be its Weil restriction along $q$, i.e. $q_*\widetilde{G}(T) := \Hom_\Xt(T \times_X \Xt, \widetilde{G})$ for every $T$ over $X$. It follows from \cite[Theorem 4 and Proposition 5, Section 7.6]{bosch1990neron} that $q_*\widetilde{G}$ is representable by a smooth group scheme over $X$. The actions of $\Gamma$ on $G$ and on $\Xt$ induce the action of $\Gamma$ on $q_*\widetilde{G}$ given by 
\[(\gamma \cdot f)(t,\xt):= \rho(\gamma)^{-1} f( \gamma(t,\xt))= \rho(\gamma)^{-1} f(t,\gamma^*(\xt))\]
for all $t \in T$ and $\xt \in \Xt$.

We define $\H$ to be the subgroup of $\Gamma$-invariants of $q_*(\widetilde{G})$, i.e. \[\mathcal{H}:=(q_*\widetilde{G})^\Gamma.\] 
We denote by $\h$ the sheaf of Lie algebras of $\H$. Since $\H$ is smooth, as shown in \cite[Proposition 3.4]{Edixhoven1992Neron}, $\h$ is a vector bundle on $X$ which is moreover equipped with a structure of Lie algebra.

\begin{eg} \label{eg-SLn} Let $\rho \colon \Gamma=\ZZ/2\ZZ \to \Aut(\SLn[r])$ be given by $\rho(\gamma)M=(M^t)^{-1}$ and $q \colon \Xt \to X$ a $\Gamma$-covering of smooth curves. The group $\H=(q_*(\SLn[r] \times \Xt))^\Gamma$ is the quasi split special unitary group associated to the extension $k(X) \subseteq k(\Xt)$. Observe that only in the case $r=2$ this action comes from inner automorphisms.
\end{eg}

\begin{rmk} The action of $\Gamma$ on $G$ via $\rho$ induces an action on $\g:=\Lie(G)$. We equivalently could have defined $\h$ as the Lie algebra of $\Gamma$-invariants of $q_*(\g \otimes_k \Oo_\Xt)$.\end{rmk}

\subsection{Properties of $\Gamma$-coverings} \label{subsec-propGammacov}

In this section we recall the definitions and properties of coverings of curves. The main refernce is \cite{bertin2007champs}, but we make the stronger assumption that all the points of $\Xt$ which are fixed by a non trivial element of $\Gamma$ are smooth.

\subsubsection{Ramification and branch divisors} Consider a $\Gamma$-covering $(f \colon \Xt \overset{q}\to X \overset{\pi}\to S)$. We define the \emph{ramification divisor} $\Rrt$ to be the effective Cartier divisor $(p-1)\Xt^\Gamma$, where $\Xt^\Gamma$ is the subscheme of $\Xt$ fixed by $\Gamma$. Equivalently, since $\Gamma$ does not have proper subgroups, $\Xt^\Gamma$ is the complement of the étale locus of $q$, which is either empty or an effective Cartier divisor of $\Xt$. The \emph{reduced branch divisor} $\Rr$ is the effective divisor given by the image of $\Xt^\Gamma$ in $X$. One can moreover observe that $q_*(\Oo(-\Rrt))^\Gamma$ is isomorphic to $\Oo(-\Rr)$.

\begin{rmk} \label{rmk-divisori} If the map $q$ is not étale both divisors $\Rrt$ and $\Rr$ are finite and étale over $S$. This is a proved in \cite[Proposition 3.1.1]{bertin2007champs} for the smooth case only and in \cite[Proposition 4.1.8]{bertin2007champs} for the general situation. \end{rmk}

The ramification divisors are naturally related to tangent sheaves of $X$ and $\Xt$. Let $\T_{\Xt/S}$ be the tangent sheaf of $\Xt$ relative to $S$, so that its sections are $f^{-1}\Oo_S$-linear derivations of $\Oo_{\Xt}$. Consider its pushforward to $X$ along $q$ and notice that the action of $\Gamma$ on $q_*\Oo_{\Xt}$ induces an action on $q_* \T_{\Xt/S}$ by sending a derivation $D$ to $\gamma D \gamma^{-1}$. The following statement, which describes the $\Gamma$-invariants of $q_* \T_{\Xt/S}$, follows from \cite[Proposition 4.1.11]{bertin2007champs}.

\begin{prop} \label{prop-tang} The sheaf $(q_*\T_{\Xt/S})^\Gamma$ over $X$ is isomorphic to $\T_{X/S}(-\Rr)$. \qed \end{prop}

\subsubsection*{Hurwitz data} The Hurwitz data provide a description of the action of $\Gamma$ at the ramification points. Before working with families we consider $q \colon \Xt \to X$, a $\Gamma$-covering of curves over $k$. Let $\xt \in \Xt(k)^\Gamma$ be a ramification point and up to the choice of a local parameter $t$ the formal disc around $\xt$ is isomorphic to $\Spec(k[\![t]\!])$. Since $\Gamma$ fixes $\xt$, one of its generators acts on $k[\![t]\!]$ by sending $t$ to $\zeta t$ for a primitive $p$-th root unity $\zeta$. It follows that the action of $\Gamma$ on $\Spec(k[\![t]\!])$ is uniquely determined by non trivial characters $\chi_{\xt} \colon \Gamma \to k^*$. Let $\Char(\Gamma)^*$ be the set of all non trivial characters of $\Gamma$ and set $R_+(\Gamma) := \oplus_{\chi \in \Char(\Gamma)^*} \ZZ \chi$. The \emph{ramification data} or \emph{Hurwitz data} of a $\Gamma$-covering $\Xt \to X$ is the element \[\xi :=\sum_{\xt \in \Xt^\Gamma}\chi_{\xt} \in R_+(\Gamma).\] 
The \emph{degree} of $\xi=\sum b_i \chi_i$ is $\deg(\xi):=\sum b_i$. Note that $\deg(\xi)=\deg(\Xt^\Gamma)=\deg(\Rr)$.

\begin{defi} Let $\Xt \to X \to S$ be a $\Gamma$-covering with $S$ connected. We say that it has Hurwitz data $\xi \in R_+(\Gamma)$ if $\xi$ is the Hurwitz data of one, hence all (see \cite[Lemme 3.1.3]{bertin2007champs}), of its fibres.
\end{defi}

We fix for the next two lemmas, a generator $\gamma$ of $\Gamma$ and $\zeta \in k$ a primitive $p$-th root of $1$. This identifies the set of characters of $\Gamma$ with $\{0, \dots, p-1\}$.

\begin{lem} \label{lem-decompE} Denote by $\E_i$ the $\Oo_{X}$-submodule of $q_*\Oo_{\Xt}$ where $\gamma$ acts by multiplication by $\zeta^i$. Then 
\[ q_*\Oo_{\Xt}= \bigoplus_{i=0}^{p-1} \E_i \AND \E_i \otimes \E_{p-i} \cong \Oo(-\Rr) \quad \text{for } i \neq 0.\]
\end{lem}

\proof The action of $\Gamma$ on $q_*\Oo_\Xt$ provides the decomposition with $\E_0 \cong \Oo_X$. Observe furthermore that since $\Oo_X$ is $\Gamma$-invariant, each eigenspace $\E_i$ is naturally an $\Oo_X$-module. Since the action of $\Gamma$ is compatible with the product in $q_*\Oo_\Xt$, the tensor product $\E_i \otimes \E_{p-i}$ is a submodule of $\E_0\cong \Oo_X$. Outside the branch divisor $\Rr$ this is an isomorphism so we only need to check what is the image along $\Rr$. Let $x \in \Rr$ and call $\xt \in \Rrt$ the point above $x$ so that $\widehat{\Oo}_{\Xt,\xt} \cong k[\![t]\!]$, with $\gamma(t)=\zeta^nt$ with $n \in \{1, \dots, p-1\}$. If follows that $(\widehat{\E}_i)_x \cong t^{[i/n]}R[\![t^p]\!]$ and $(\widehat{\E}_{p-i})_x \cong t^{p-[i/n]}k[\![t^p]\!]$, where $[i/n] \in \{1, \dots, p-1\}$ denotes the product of $i$ with the multiplicative inverse of $n$ in $\ZZ/p\ZZ$. It follows that $(\widehat{\E}_i)_x \otimes (\widehat{\E}_{p-i})_x \cong t^pk[\![t^p]\!]$ which is isomorphic to the completion of $\Oo(-\Rr)$ at the point $x$.
\endproof

\begin{lem} \label{lem-decomp} Denote by $\g^{(i)}$ the submodule of $\g$ where $\gamma$ acts by multiplication by $\zeta^i$. The sheaf $\h$ decomposes as \[ \h= \bigoplus_{i=0}^{p-1}\g^{(-i)} \otimes_k \E_i.
\]
\end{lem}

\proof As the action of $\gamma$ on $\g$ is diagonalizable with eigenvalues belonging to $\{1,\zeta, \cdots, \zeta^{p-1}\}$, we can decompose $\g$ as $\bigoplus \g^{(-i)}$. As $\h$ is the Lie algebra of $\Gamma$-invariants of $q_*(\Oo_{\Xt} \otimes_k \g)= q_*\Oo_{\Xt} \otimes_k \g$, we can combine this with the description of $q_*\Oo_{\Xt}$ provided by Lemma \ref{lem-decompE} to obtain the wanted decomposition of $\h$.
\endproof

\begin{rmk} In view of the conclusions of Remark \ref{rmk-loccoord}, the above lemmas, which we have proved in the case of coverings of curves over $\Spec(k)$, hold for any family of coverings or curves over an arbitrary scheme $S$.\end{rmk}

\begin{eg} The previous Lemma shows in particular that when the group $\Gamma$ acts trivially on $G$, then $\h$ is isomorphic to $\g\otimes \Oo_X$. It follows that when we are in this situation we retrieve the classical construction of conformal blocks attached to simple Lie algebras. 
\end{eg}

\subsection{Hurwitz stacks} \label{subsec-Hurwitz}

We define in this section the stack parametrizing $\Gamma$-coverings with fixed Hurwitz data $\xi \in R_+(\Gamma)$.
Let $g$ be a non negative integer. Let $f \colon \Xt \overset{q}\to X \overset{\pi}\to S$ be a $\Gamma$-covering of curves and let $\sigma \colon S \to X$ be a section of $\pi$ with $\sigma(S)$ disjoint from the nodes of $X$ and from the branch locus $\Rr$ of $q$. We say that the covering is \emph{stably marked} by $\sigma$ if $(X, \sigma \cup \Rr)$ is a stable marked curve \cite[Définition 4.3.4. and Proposition 5.1.3]{bertin2007champs}. This means that $X$ is a family of curves with at most nodal singularities, the section $\sigma$ is disjoint from the nodes and from the branch locus $\Rr$, and the automorphism group of each fiber of $X \to S$ preserving the marked locus $\sigma \cup \Rr$ is finite. 

\begin{defi} We define the \emph{Hurwitz stack} $\bHur[n]$ as
\[ \bHur[n](S)=\left\langle f \colon \Xt \overset{q}\to X \overset{\pi}\to S, \:\{\sigma_j \colon S \to X\}_{j=1}^n \text{ such that i and ii hold}\right\rangle
\]
\begin{enumerate}[label=\roman*.]
\item the map $q \colon \Xt \to X$ is a $\Gamma$-covering of curves with ramification data $\xi$;
\item $(X,\{\sigma_j\})$ is an $n$-marked curve of genus $g$ such that each $\sigma_j(S)$ is disjoint from the nodal locus and from the branch divisor $\Rr$ and such that the covering is stably marked by $\{\sigma_j\}$.
\end{enumerate}\end{defi}

When $n=0$ we omit the subscript and use the notation $\bHurz$. We denote by $\Hur[n]$ the open substack of $\bHur[n]$ parametrizing $\Gamma$-coverings of smooth curves.

\begin{rmk} Although the notation might suggest that $\bHur[n]$ is a compactification of $\Hur[n]$, the stack $\bHur[n]$ is not proper because the ramification points avoid both marked and singular points. \end{rmk}

In the previous section we explained how to associate to each $\Gamma$-covering $(\Xt \overset{q}\to X \overset{\pi}\to S) \in \bHurz(S)$, a group $\H$ (resp. a sheaf of Lie algebras $\h$) over $X$. This defines a group $\H_{\univ}$ (resp. a sheaf of Lie algebras $\h_{\univ}$) on $\Cun$, where we denote by $\Ctun \to \Cun$ the universal covering on $\bHurz$. The same construction works on $\bHur[n]$, defining then $\H_{\univ}$ and $\h_{\univ}$ on the universal curve $\Cun$ of $\bHur[n]$.

\begin{rmk} \label{rmk-NCD} The complement $\Delta_{\univ}:= \bHur[n] \setminus \Hur[n]$ is a normal crossing divisor. First of all observe that $\Delta_{\bMgn[d]} := \bMgn[d] \setminus \Mgn[d]$ is a normal crossing divisor: in fact given a nodal curve $X \to \Spec(k)$ with a reduced divisor $D$ of degree $d$, there exists a versal deformation $\Xx \to S$ where the locus $\Delta \subset S$ consisting of singular curves is a normal crossing divisor of $S$ \cite{ArbarelloCornalba}. We now want to compare the deformation theory of a $\Gamma$-covering $(\Xt \to X, \{\sigma_i\})$ to the one of $(X, \{\sigma_i\} \cup \Rr)$. Following \cite[Théorème 5.1.5]{bertin2007champs} we see that the natural map $\delta \colon \Def(\Xt \to X, \{\sigma_i\}) \to \Def(X, \{\sigma_i\} \cup \Rr)$ fails to be an isomorphism only when the intersection between $\Rr$ and $X^{sing}$ is not empty, but since by assumption we impose that $\Rr \cap X^{sing}= \emptyset$, in our contest this map is always an isomorphism. This then allow to obtain, from the versal deformation $\Xx \to S$ of $(X, \{\sigma_i\} \cup \Rr)$, the versal deformation $(\Xxt \to \Xx, \{\varsigma_i\})$ of $(\Xt \to X, \{ \sigma_i\})$, and hence deduce from the theory of $\bMgn[n+\deg(\Rr)]$ that $\Delta_{\univ}$ is a normal crossing divisor.
\end{rmk}

The following statement, which is given by \cite[Proposition 2.3.9. and Théorème 6.3.1]{bertin2007champs}, describes the properties of the above stacks. 

\begin{prop} The stacks $\bHur[n]$ and $\Hur[n]$ are connected smooth Deligne-Mumford stacks of finite type over $\Spec(k)$.\end{prop}

\begin{rmk} We want to remark that the role of the ramification data is to guarantee the connectedness of $\bHur[n]$ and $\Hur[n]$ \cite[Proposition 2.3.9]{bertin2007champs}. If the group $\Gamma$ were not cyclic, however, fixing the ramification data would in general not suffice to guarantee the connectedness of $\bHur[n]$ or $\Hur[n]$.\end{rmk}

Instead of marking the curve $X$, we can mark the curve $\Xt$, so that we define.

\begin{defi} For each $k$-scheme $S$ we set
\[ \bHurs{n}{\xi}(S)=\left\langle f \colon \Xt \overset{q}\to X \overset{\pi}\to S, \:\{\tau_j \colon S \to \Xt \}_{j=1}^n \text{ such that i and ii hold}\right\rangle
\]
\begin{enumerate}[label=\roman*.]
\item the map $q \colon \Xt \to X$ is a $\Gamma$-covering of curves with ramification data $\xi$;
\item $(\Xt,\{\tau_j\})$ is an $n$-marked curve with $q \tau_j(S)$ pairwise disjoint, $q \tau_j(S)$ disjoint from $\Rr$ for all $j$ and such that the covering $q$ is stably marked by $\{q\tau_j\}$.
\end{enumerate}
\end{defi}

It follows, from the fact that the image of $\tau$ lies in the étale locus of $q$, that the map \[\textbf{Forg}^n_n \colon \Hurs{n}{\xi} \to \Hur[n], \qquad (\Xt \overset{q}\to X \overset{\pi}\to S, \:\{\tau_j\}) \mapsto (\Xt \overset{q}\to X \overset{\pi}\to S, \: \{q\tau_j\})
\] is an étale and surjective morphism of stacks. For any $n,m \in \ZZ_{\geq 0}$ we also have the forgetful map \[\textbf{Forg}_{n+m,n} \colon\Hur[n+m] \to \Hur[n].\]

The advantage of marking the curve $\Xt$ instead of $X$ lies in the following proposition. 

\begin{prop} \label{TauIso} Let $(\Xt \overset{q}\to X \overset{\pi}\to S, \tau) \in \Hurs{1}{\xi}(S)$ and write $\sigma=q \tau$. The section $\tau$ induces an isomorphism between ${\sigma}^*\H$ and $G \times_k S$.
\end{prop}

\proof Construct the cartesian diagram \[ \xymatrix{\widetilde{S} \ar[r]^{\widetilde{\sigma}} \ar[d]^{q_S}& \Xt \ar[d]^q\\
S \ar[r]^{\sigma} & X}\]
and since by assumption the image of $\sigma$ lies in the étale locus of $q$ the left vertical arrow $q_S$ is étale and it has a section given by $\tau$. This implies that $\widetilde{S}$ is isomorphic to $\coprod_{\gamma_i \in \Gamma} S$.  Observe that ${q_S}_*\widetilde{\sigma}^*(\widetilde{G})\cong \sigma^*q_*(\widetilde{G})$ and that taking $\Gamma$-invariants commutes with restriction along $\sigma$. It follows that
\[\sigma^*\mathcal{H}=\left(\sigma^*q_*(\widetilde{G})\right)^\Gamma =\left({q_S}_*\widetilde{\sigma}^*(\widetilde{G})\right)^\Gamma =\left( {q_{S}}_* \left(\coprod_{\gamma_i \in \Gamma} S\times G\right) \right)^\Gamma  =  \left(\prod_{\gamma_i \in \Gamma} S \times G \right)^\Gamma\] 
where $\gamma_j \in \Gamma$ acts on $\prod_{\gamma_i \in \Gamma} S\times G$ by sending $(s_i,g_i)_{\gamma_i}$ to $(s_i,\gamma_j(g_i))_{\gamma_j\gamma_i}$. It follows that the invariant elements are of the form $(s, \gamma_i(g))_{\gamma_i}$ for any $s\in S$ and $g \in G$, so that the projection on any component of $S\times G$ realized an isomorphism between $\sigma^*\mathcal{H}$ and $G \times S$. The map $\tau$ selects a preferred component, giving in this way a canonical isomorphism. 
\endproof

\section{The sheaf of covacua and of conformal blocks} \label{sec-ConformalBlocks}

In this section we give the definition of the sheaves of covacua and of conformal blocks on $\bHur[n]$. We begin by considering the case $n=1$ and explain how to construct the sheaf $\VV_\ell(\V)$ on $\bHur$ attached to a representation $\V$ of $\sigma^*\h_{\univ}$. The sheaf of conformal blocks will be realized as the dual of $\VV_\ell(\V)$. In order to define this sheaf on $\bHur$, we will define it for any family $(f \colon \Xt \overset{q}\to X \overset{\pi}\to S, \sigma)$ over an affine and smooth base $S=\Spec(R)$. We can assume moreover that $X \setminus \sigma(S) \to S$ is affine and we will see in Remark \ref{rmk-dropaffine} how to drop this assumption.

For the classical definition of the sheaf of conformal blocks attached to a representation of $\g$ one can refer to \cite{tsuchiya1989conformal} or to \cite{looijenga2005conformal}. We will use the latter as main reference, from which we borrow the notatation.

Let $X^*:=X \setminus \sigma(S)$ and denote by $\Aa$ the pushforward to $S$ of $\Oo_{X^*}$, i.e. \[\Aa:=\pi_*{j_A}_*\Oo_{X^*}\] where $j_A$ denotes the open immersion $X^* \to X$. Since the map $\pi$ restricted to $X^*$ is affine we have that $X^*=\SPec(\Aa)$ and that $\Aa =\pi_* \varinjlim_n \mathcal{I}_\sigma^{-n} = \varinjlim_n \pi_* \mathcal{I}_\sigma^{-n}$ where $\mathcal{I}_\sigma=\Oo_X(-\sigma(S))$ is the ideal defining $\sigma(S)$.

We denote by $\Op$ the formal completion of $\Oo_X$ along $\sigma(S)$: by definition $\sigma$ gives a short exact sequence \[ 0 \to \mathcal{I}_\sigma \to \Oo_X \to \Oo_{\sigma(S)} \to 0 \] of $\Oo_X$-modules. We define \[\Op:= \pi_*\varprojlim_{n \in \NN} \Oo_X \Qdx (\mathcal{I}_\sigma)^n = \varprojlim_{n \in \NN} \pi_*\Oo_X \Qdx (\mathcal{I}_\sigma)^n\] which is naturally a sheaf of $\Oo_S$-modules. We denote by $\Ll$ the $\Oo_S$-module \[\Ll:= \varinjlim_{N \in \ZZ_{\geq 0}} \pi_* \varprojlim_{n\in \ZZ_{> 0}} \mathcal{I}_\sigma^{-N} / \mathcal{I}_\sigma^{n}\] which is equipped with a natural filtration \[F^N\!\Ll=\pi_* \varprojlim_{n \in \ZZ_{> 0}} \mathcal{I}_\sigma^{N} / \mathcal{I}_\sigma^{N+n} \qquad \text{ for }N \geq 0\] and  \[F^N\!\Ll=\pi_* \varprojlim_{n \in \ZZ_{> 0}} \mathcal{I}_\sigma^{N} / \mathcal{I}_\sigma^{n}\qquad \text{ for }N \leq -1\] taking into account the order of the poles or zeros along $\sigma(S)$. 

\begin{rmk} \label{rmk-loccoord} Recall that when $R=k$, the choice of a local parameter $t$, i.e. of a generator of $\mathcal{I}_\sigma$, gives an isomorphism $\Op \cong k[\![t]\!]$ and hence $\Ll\cong k(\!(t)\!)$ and so $F^n\!\Ll \cong t^n k[\![t]\!]$. In the general case, since $\mathcal{I}_\sigma$ is locally principal, for every $s \in \sigma(S)$ we can find an open covering $U$ of $X$ containing $s$ and such that $\mathcal{I}_\sigma|_{U}$ is principal. Let denote by $S'$ the open of $S$ given by $\sigma^{-1}(U)$ and by $U'$ the open $U \cap \pi^{-1}{S'}$. 
Then $\mathcal{I}_\sigma|_{U'}$ is principal and $\varprojlim_n \Oo_{U'} / (\mathcal{I}_\sigma|_{U'})^n$ is isomorphic to $\Oo_{S'}[\![t]\!]$, where $t$ is a generator of $\mathcal{I}_\sigma|_{U'}$. This moreover implies that the completion of $\Op$ at a point $s\in S$ is isomorphic to $\widehat{\Oo}_{S,s}[\![t]\!]$, where $\widehat{\Oo}_{S,s}$ denotes the completion of $\Oo_S$ at $s$.\end{rmk}

Denote by $\hA$ the restriction of $\h$ to the open curve $X^*$ and by $\hL$ the "restriction of $\h$ to the punctured formal neighbourhood around $\sigma(S)$", and consider both sheaves as $\Oo_S$-modules naturally equipped with a Lie bracket. In other words we set 
\begin{align*} \hA & :=\pi_*{j_A}_* {j_A}^*(\h)=  \pi_* \left(\varinjlim_{N \in \ZZ_{\geq 0}} \mathcal{I}_\sigma^{-N} \otimes_{\Oo_X} \h \right)= \varinjlim_{N \in \ZZ_{\geq 0}} \pi_* ( \mathcal{I}_\sigma^{-N} \otimes_{\Oo_X} \h)\\
\hL &:= \varinjlim_{N \in \ZZ_{\geq 0}} \pi_* \varprojlim_{n \in \ZZ_{> 0}} \mathcal{I}_\sigma^{-N} / \mathcal{I}_\sigma^{n} \otimes_{\Oo_X} \h. 
\end{align*}  

The following observations follow from the definitions.\begin{enumerate}
\item The injective morphism $\mathcal{I}_\sigma^{-N} \to \varprojlim_n \mathcal{I}_\sigma^{-N} / \mathcal{I}_\sigma^n$ induces the inclusion $\hA \to \hL$. 
\item The filtration on $\Ll$ defines the filtration $F^*\hL$ as
\[F^N(\hL)= \pi_* \varprojlim_{n \in \ZZ_{> 0}} \mathcal{I}_\sigma^{N} / \mathcal{I}_\sigma^{N+n} \otimes_{\Oo_X} \h \AND F^{-N}(\hL)= \pi_* \varprojlim_{n \in \ZZ_{> 0}} \mathcal{I}_\sigma^{-N} / \mathcal{I}_\sigma^{n} \otimes_{\Oo_X} \h
\]for all $N \in \ZZ_{\geq 0}$ and we denote $F^0(\hL)$ by $\h_{\Op}$.
\item We could have equivalently defined $\hA$ as the Lie subalgebra of $\Gamma$-invariants of $f_*(\g \otimes_k {j_{\widetilde{A}}}_* \Oo_{\Xt*})$ where $j_{\widetilde{A}}$ denotes the open immersion of $\Xt^*:=\Xt \times_X X^* \to \Xt$. This follows from the equalities
\[{j_A}^* \h= {j_A}^*(q_*(\g \otimes_k \Oo_{\Xt}))^\Gamma= ({j_A}^*q_*(\g \otimes_k \Oo_{\Xt}))^\Gamma = q_* ({j_{\widetilde{A}}}^*(\g \otimes_k \Oo_X^*))^\Gamma.
\]
Similarly $\hL$ is the Lie subalgebra of $\Gamma$-invariants of $\g \otimes_k \widehat{\Ll}$, where \[\widehat{\Ll}:=\varinjlim_{N \in \ZZ_{\geq 0}} f_* \varprojlim_{n \in \ZZ_{> 0}} (\g \otimes_k q^*(\mathcal{I}_\sigma^{-N}) / q^*(\mathcal{I}_\sigma^n)).\]
\end{enumerate}

\begin{rmk} \label{rmk-IsomhLgL} Since $\sigma(S)$ has trivial intersection with $\Rr$, we can find an étale cover of $S$ such that $q^{-1}(\sigma(S))= \coprod_{\Gamma} S$ or in other terms such that the pull back of $\mathcal{I}_\sigma$ to the cover totally splits, i.e. $q^*\mathcal{I}_\sigma= \prod_{\gamma\in \Gamma} \mathcal{I}_{\sigma,\gamma}$. This implies that
\[\hL \cong  \left( \g \otimes_k \bigoplus_{\gamma \in \Gamma} (\varinjlim_N f_* \varprojlim_n \mathcal{I}_{\sigma,\gamma}^{-N} / \mathcal{I}_{\sigma,\gamma}^n ) \right)^\Gamma\]
which leads to $\hL \cong (\g \otimes_k (\oplus_{\gamma \in \Gamma}\Ll))^\Gamma$ where the action is given by
\[\gamma_j * \left((X_{\gamma}f_{\gamma})_{\gamma}\right) = (\gamma_j(X_{\gamma})f_\gamma)_{\gamma_j \gamma} \quad \text{for all } X_\gamma \in \g \text{ and } f_\gamma \in \Ll.\]
It follows that the invariant elements are combination of elements of the type $(\gamma(X)f)_{\gamma}$ for $X \in \g$ and $f \in \Ll$. For every $\gamma \in \Gamma$, the projection on the $\gamma$-th component
\[\text{pr}_\gamma \colon \hL \to \gL:=\g \otimes_k \Ll, \quad  (\gamma(X)f)_{\gamma} \mapsto \gamma(X)f\]
defines a non canonical isomorphism of sheaves of Lie algebras of $\hL$ with $\gL$. The inverse of $pr_\gamma$ is the map that sends the element $Xf$ of $\gL$ to the $p$-tuple $(\gamma_j(\gamma^{-1}(X))f)_{\gamma_j}$.
\end{rmk}

\subsection{The central extension of $\hL$} \label{subsec-CentralExtensionhL}

Once we have defined $\hL$ and $\hA$, in order to define $\VV_\ell(\V)$, we need to extend $\hL$ centrally. Following \cite[Chapter 7]{kac1994infinite}, \cite{tsuchiya1989conformal} and \cite{looijenga2005conformal} we construct this central extension using a normalized Killing form and the residue pairing.

\subsubsection*{Normalized Killing form}

We fix once and for all a maximal torus $T$ of $G$ and a Borel subgroup $B$ of $G$ containing $T$, or equivalently we fix the root system $R(G,T)=R(\g,\mathfrak{t}) \subseteq \mathfrak{t}^\vee:=\Hom(\mathfrak{t},k)$ of $G$ and a basis $\Delta$ of positive simple roots, where $\mathfrak{t}=\text{Lie}(T)$. Given a root $\alpha$ we denote by $H_\alpha \in \mathfrak{t}$ the associated coroot. 

Denote by $(\,|\,)\colon \g \otimes \g \to k$ the unique multiple of the Killing form such that $(H_\theta|H_\theta)=2$ where $\theta$ is the highest root of $\g$. As $\g$ is simple, this form gives an isomorphism $(\,|\,)$ between $\g$ and $\g^\vee:=\Hom(\g, k)$. 
Pulling back this form to $\Xt$, then pushing it forward along $q$ we obtain 
\[ q_*\tilde{(\,|\,)}\colon q_*(\widetilde{\g}) \otimes q_*(\widetilde{\g}) \to q_*(\Oo_\Xt)
\]
where $\widetilde{\g}:=\g \otimes \Oo_\Xt$. Since the Killing form is invariant under automorphisms of $\g$, the bilinear form $q_*\tilde{(\,|\,)}$ is $\Gamma$-equivariant. Taking $\Gamma$-invariants we obtain the pairing
\[(\,|\,)_\h \colon \h \otimes_{\Oo_X} \h \to \Oo_X
\]
which however is not perfect because of ramification. Combining this with the multiplication map $\mathcal{I}_\sigma^{-N}/\mathcal{I}_\sigma^{N+n} \times \mathcal{I}_\sigma^{-N}/\mathcal{I}_\sigma^{N+n} \to \mathcal{I}_\sigma^{-2N}/\mathcal{I}_\sigma^n$ and taking the limit on $n$ and $N$ we obtain the pairing $(\,|\,)_\hL \colon \hL \otimes_{\Ll} \hL \to \Ll$ which is perfect. 

\subsubsection*{Residue pairing} \label{subsubsec-ResiduePairing}

We introduce the sheaf $\theta_{\Ll/S}$ of continuous derivations of $\Ll$ which are $\Oo_S$ linear. Denote its $\Ll$-dual by $\omega_{\Ll/S}$: this is the sheaf of continuous differentials of $\Ll$ relative to $\Oo_S$. Observe that when $\Op \cong R[\![t]\!]$ we have that $\theta_{\Ll/S}$ is isomorphic to $R(\!(t)\!) d/dt$ and $\omega_{\Ll/S}$ to $R(\!(t)\!)dt$.

The residue map $\Res \colon \omega_{\Ll /S} \to \Oo_S$ is computed locally as $\Res( \sum_{i\geq N} \alpha_it^idt)=\alpha_{-1}$. Composing this with the canonical morphism $\hL^\vee \times \hL \to \Ll$ we obtain the perfect pairing
\[\Res_\h \colon \omega_{\Ll/S}\otimes_\Ll \hL^\vee \times \hL \to \Oo_S.
\] 

\subsubsection*{The differential of a section}
Let $d \colon \Oo_\Xt \to \Omega_{\Xt / S}$ be the universal derivation, which induces the morphism $d \colon \g \otimes_k \Oo_{\Xt} \to   \Omega_{\Xt / S}\otimes_k\g$ by tensoring it with $\g$. Let $U \subseteq X \setminus \Rr$ be an open subscheme of $X$ containing $\sigma(S)$ and which is smooth over $S$, and call $\widetilde{U}=U \times_X \Xt$. Once we restrict $d$ to $\widetilde{U}$ and we push it forward along $q$ we obtain the map \[d \colon q_*(\g \otimes_k \Oo_{\widetilde{U}}) \to \Omega_{U/S} \otimes_{\Oo_U}q_*(\g \otimes_k \Oo_{\widetilde{U}})  \] by using the projection formula. Taking $\Gamma$-invariants one obtains $d \colon \h|_{U} \to   \Omega_{U/S}\otimes_{\Oo_U} \h|_{U}$ and since $\sigma(S) \subset U$, this induces the map $d \colon \hL \to \omega_{\Ll/S} \otimes_\Ll \hL$. We can furthermore compose this map with the morphism $\hL \to \hL^\vee$ given by the normalized Killing form $(\,|\,)_\hL$, obtaining \[d_\hL \colon \hL \to \omega_{\Ll/S} \otimes_\Ll \hL^\vee.\]

\begin{rmk} We could have equivalently defined $d_\hL$ by using the local isomorphism between $\hL$ and $\gL$. Using this approach, we can describe $d_\hL$ as the map which associates to the element $Xf \in \gL$, the element $df \otimes (X|-)$ belonging to $\omega_{\Ll/S} \otimes_\Ll (\gL)^\vee$.
\end{rmk}

\begin{rmk} \label{rmk-dkill} Given $X,Y \in \hL$, we simply write $(dX|Y)$ for $d_\hL(X)(Y) \in \omega_{L/S}$. Note that the following equality holds $d_{\Ll}(X|Y)_\hL= (dX|Y) + (X|dY)$, where $d_{\Ll} \colon \Ll \to \omega_{\Ll/S}$ is the universal derivation.\end{rmk}

\subsubsection*{The central extension of $\hL$}

We have introduced all the ingredients we needed to be able to define the central extension $0 \to c\Oo_S \to \HhL \to \hL \to 0$ of $\hL$ where $c$ is a formal variable.

\begin{defi} We define the sheaf of Lie algebras $\HhL$ to be $\hL \oplus c\Oo_S$ as $\Oo_S$-module, with $c\Oo_S$ being in the centre of $\HhL$ and with Lie bracket defined as
\[ [X,Y]:=[X,Y]_{\hL} + c \Res_\h\left(d_{\hL}(X) \otimes Y\right)=[X,Y]_{\hL} + c \Res\left(dX|Y\right)
\]for all $X,Y \in \hL$.
\end{defi}

The Lie algebra $\HhL$ comes equipped with the filtration $F^i\HhL=F^i\hL$ for all $i \geq 1$ and $F^i\HhL=F^i\hL \oplus c\Oo_S$ for $i \leq 0$.

\begin{rmk} \label{rmk-HhLandHgL} We can locally describe the Lie algebra $\HhL$ as follows. Locally on $S$ we know that we can lift $\sigma$ to $p$ distinct sections $\tau_0, \dots, \tau_{p-1}$ of $f$ so that $\widetilde{\Ll}$ is isomorphic to $\bigoplus_{i=0}^{p-1} \Ll_i$. We can define the central extension $\widehat{\g}\widetilde{\Ll}$ of $\g \otimes_k \widetilde{\Ll}$ as the module $\g \otimes_k \widetilde{\Ll} \oplus \widetilde{c} \Oo_S$ with the following Lie bracket:
\[ [ (X_i f_i)_i , (Y_i g_i)_i] = ([X_i,Y_i]f_ig_i)_i \oplus \widetilde{c} \sum_{i=0}^{p-1} (X_i|Y_i)\Res_i(df_i g_i) 
\]
where $\Res_i$ computes the residue at the section $\tau_i$. The induced action of $\Gamma$, which acts trivially on $\widetilde{c}$, respects the Lie bracket, hence the $\Gamma$-invariants define a central extension of $\hL$ which coincide with $\HhL$ by setting $c=\widetilde{c}/p$.
\end{rmk}

As $\hA \subset \hL$, one might wonder which is the Lie algebra structure induced on $\widehat{\h}_\Aa$. 

\begin{prop} The inclusion $\hA \subset \hL$ induces a natural inclusion of $\hA$ in $\HhL$.
\end{prop} 

\proof The map $\hA \to \HhL$ is given by the inclusion of $\hA$ into $\hL$ as modules, so we are only left to prove that this is a Lie algebra morphism. This can be checked locally on $S$, so in view of the Remark \ref{rmk-HhLandHgL} $\HhL$ is given by the $\Gamma$-invariants of $\widehat{\g}\otimes \widetilde{\Ll}$. By definition $\hA$ is given by the $\Gamma$-invariants of $\g \otimes \widetilde{A}$ and thanks to \cite[Lemma 5.1]{looijenga2005conformal} we know that $\g \otimes \widetilde{A}$ is a Lie sub algebra of $\widehat{\g}\otimes \widetilde{\Ll}$.
\endproof

We can also prove the previous proposition as a consequence of the following two lemmas, which we present as they are going to be useful in Section \ref{sec-ProjConn}. Let denote by $\omega_{\Aa/S}$ the pushforward to $S$ of $\omega_{X\setminus \sigma(S)/S}$, the relative dualizing sheaf of $X\setminus \sigma(S)$ over $S$. This is a subsheaf of $\omega_{\Ll/S}$.

\begin{lem} The image of $\hA$ via $d_{\hL}$ is $\omega_{\Aa/S} \otimes_\Aa \hA^\vee$.
\end{lem}

\proof We can restrict to the case of family of smooth curves, as on the singular points the result follows from \cite[Lemma 5.1]{looijenga2005conformal} by identifying $\h$ with $\g \otimes \Oo_X$. Recall from Lemma \ref{lem-decomp} that $\h=\bigoplus_{i=0}^{p-1} \g^{(-i)}\otimes_k \E_i$, and note that the image of $\E_i$ under $d$ is $\E_{i}(\Rr)\otimes_{\Oo_X} \Omega_{X/S}$ for $i \neq 0$. To check this fact, it is enough to consider what happens locally at a point $x \in \Rr$. As in the proof of Lemma \ref{lem-decompE}, the completion of $\E_i$ at $x$ is isomorphic to $t^{[i/n]}R[\![ t^p]\!]$, and the image of $t^{[i/n]}$ under $d$ is $[i/n] t^{[i/n]-1}dt$ and $dt= p^{-1} t^{1-p}d(t^p)$, we conclude that the image is a scalar multiple of $t^{[i/n]-p}d(t^p)$, hence it belongs to the completion of $\E_i(\Rr)\otimes \Omega_X$ at $x$. Observe furthermore that $(\,|\,)$ gives an isomorphism between $\g^{(i)}$ and the dual of $\g^{(-i)}$. Since $\E_i \otimes_{\Oo_X} \E_{p-i} \cong \Oo(-\Rr)$ for $i \neq 0$ and $\E_0= \Oo_X$ it follows that if $i \neq 0$ we have
\begin{align*} d_\hL (\g^{(-i)}\otimes_k \E_i) &= ( \g^{(-i)} \otimes_k \E_i(\Rr),| -)_{\h} \otimes_{\Oo_X} \Omega_{X/S} \\ &= (\g^{(-i)}| -) \otimes_k (\E_{i}(\Rr))^\vee \otimes_{\Oo_X} \Omega_{X/S} \\ &= \g^{(i)} \otimes_k \E_{-i} \otimes_{\Oo_X} \Omega_{X/S} \end{align*}
and similarly $d_\hL (\g^\Gamma \otimes_k \Oo_X) = \g^\Gamma \otimes_k \Omega_{X/S}$ which together yield $d_\hL(\hA)=\omega_{\Aa/S} \otimes_\Aa \hA^\vee$.
\endproof

\begin{lem} \label{lem-annullatore} The annihilator of $\hA$ with respect to the pairing $\Res_\h$, denoted $\text{Ann}_{\Res_\h}(\hA)$, is $\omega_{\Aa/S} \otimes \hA^\vee$.
\end{lem}

\proof Before starting with the proof, we remark that this lemma holds if we replace $\hA$ with any vector bundle $\E$ on $X$ as it is essentially a consequence of Serre duality. We start by giving a description of the quotient $\hA \Qsx \hL$, as the annihilator of $\hA$ will be the dual of that quotient with respect to the residue pairing. The double quotient $\hA \Qsx \hL \Qdx F^n\hL$ computes $R^1\pi_*(\h \otimes_{\Oo_X} \mathcal{I}_\sigma^n)$. It follows that the projective limit $\varprojlim_{n \geq 1} R^1\pi_*(\h \otimes_{\Oo_X} \mathcal{I}_\sigma^n))$ equals $\varprojlim_{n \geq 1} \hA \Qsx \hL \Qdx F^n\hL$ which is $\hA \Qsx \hL$.
As the residue pairing gives rise to Serre duality, we know that $R^1\pi_*(\h \otimes_{\Oo_X} \mathcal{I}_\sigma^n)$ is isomorphic to the dual of $\pi_*(\omega_{X/S} \otimes (\h\otimes \mathcal{I}_\sigma^n)^\vee)$. It follows that \[\text{Ann}_{\Res_\h}(\hA) = \varinjlim_{n\geq 1} \pi_*\left(\omega_{X/S} \otimes_{\Oo_X} (\h\otimes_{\Oo_X} \mathcal{I}_\sigma^n)^\vee\right).\]
which equals $\omega_{\Aa/S} \otimes_{\Aa} \hA^\vee$. \endproof

\subsection{Conformal blocks attached to integrable representations} \label{subsec-defHhn}

We have all the ingredients to define the sheaf of conformal blocks. Let $U\HhL$ denote the universal enveloping algebra of $\HhL$ and recall that $F^0\hL= \pi_* \varprojlim_n \Oo_X/ \mathcal{I}_\sigma^n$, i.e. it's the subalgebra of $\hL$ which has no poles along $\sigma(S)$. Observe that this implies that it is also a Lie sub algebra of $\HhL$.

\begin{defi} \label{def-Verma} For any $\ell \in \ZZ_{> 0}$ we define the \emph{Verma module of level $\ell$} to be the left $U\HhL$-module given by
\[\tHhe  := U\HhL \Qdx \left(U\HhL \circ F^0\hL, \, c=\ell \right).\]
\end{defi}

For what follows we will need a generalization of this module attached to certain representations of $\sigma^*\h$.

\begin{defi} An irreducible finite dimensional representation $\mathcal{V}$ of $\sigma^*\h$ is a locally free $\Oo_S$-module which is equipped with an action of $\sigma^*\h$ which locally étale on $S$, and up to an isomorphism of $\sigma^*\h$ with $\g \times \Oo_{S}$, is isomorphic to $V \otimes_k \Oo_{S}$ for an irreducible finite dimensional representation $V$ of $\g$.\end{defi}

Let $\mathcal{V}$ be an irreducible finite dimensional representation of the Lie algebra $\sigma^*\h$: we will see how this induces a representation of $\HhL$ with the central element acting as multiplication by $\ell \in \ZZ_{> 0}$. As first step, note that the exact sequence 
\[0 \to \mathcal{I}_\sigma \to \Oo_X \to \Oo_S \to 0\]
defining $\sigma(S)$ gives rise to the map of Lie algebras $[\star]_{\mathcal{I}_\sigma} \colon F^0\hL \to \sigma^*\h$ which is induced by the truncation map $\varprojlim_{n \in \ZZ_{> 0}} \Oo_X /\mathcal{I}_\sigma^n \to \Oo_X/\mathcal{I}_\sigma$.
The action of $\sigma^*\h$ on $\mathcal{V}$ is then extended to the action of $F^0\HhL =F^0\hL \oplus c \Oo_S$ by imposing, for every $v \in \mathcal{V}$ and for every $X \in F^0\hL$, the relations
\[ c * v:=\ell v \AND X * v:= [X]_{\mathcal{I}_\sigma} v.\]

In view of this, once we fix $\ell \in \ZZ_{> 0}$ we always view a representation $\V$ of $\sigma^*\h$ as a $UF^0\HhL$-module with the central part acting by multiplication by $\ell$.

\begin{defi} \label{def-VermaV} For every $\ell \in \ZZ_{> 0}$ we define the \emph{Verma module of level $\ell$ attached to $\V$} to be left $U\HhL$-module of level $\ell$ attached to $\mathcal{V}$, meaning
\[ \tHh{\mathcal{V}} := U\HhL \underset{UF^0\HhL}{\otimes} \mathcal{V}
\]where $F^0\HhL$ acts on $U\HhL$ by multiplication on the right and $U\HhL$ acts on $\tHh{\mathcal{V}}$ by left multiplication.
\end{defi}

\begin{rmk} Note that when $\V$ is the trivial representation of $\sigma^*\h$, we obtain that $\tHh{\V}$ coincides with $\tHhe$ given in Definition \ref{def-Verma}.
\end{rmk}

In the constant case $\sigma^*\h\cong \g$, the properties of $\tHh{\V}$ have been studied in \cite[Chapter 7]{kac1994infinite} when $R=k$ and $\V$ an irreducible representation of $\g$ of level at most $\ell$, where it is shown that it has a maximal irreducible quotient $\Hhn{\V}$. From this, one generalizes the construction to families of curves, but still in the constant case $\sigma^*\h\cong \g$, meaning working on $\bHurs{1}{\xi}$. The new step is to descend $\Hhn{\V}$ from $\bHurs{1}{\xi}$ to $\bHur$. 

We first of all recall the construction in the constant case in Section \ref{subsubsec-IntegrReprLevEllHurs} and then show how it descends to $\bHur$ in Section \ref{subsubsec-IntegrReprLevEllHur}.

\subsubsection{Integrable representations of level $\ell$ on $\bHurs{1}{\xi}$} \label{subsubsec-IntegrReprLevEllHurs}

The morphism $\textbf{Forg}^1_1 \colon \bHurs{1}{\xi} \to \Hur$ is a finite étale covering, so if we want to define a-module on $\bHur$, we could first define it on $\bHurs{1}{\xi}$ and later show that the construction is $\Gamma$-equivariant, hence it descends to a module on $\bHur$. As already explained in Proposition \ref{TauIso}, the advantage of working on $\bHurs{1}{\xi}$, is the identification of $\hL$ with $\gL$, which allows us to use representation theory of $\g$ and of the affine Lie algebra $\HgL$ \cite[Chapter 7]{kac1994infinite}.

We recall here some facts about representation theory of $\g$ and $\HgL$. Let $R(G,T)=R(\g,\mathfrak{t})$ be the root system of $\g$ with basis of positive roots $\Delta$. Denote by $P_+$ the set of dominant weights of $\g$ and $H_\theta$ the highest coroot of $\g$. Then for every $\ell \in \ZZ_{> 0}$ we set  
\[P_\ell := \{\lambda \in P_+ | \lambda(H_\theta)\leq \ell \}.
\]

In view of the correspondence between weights and irreducible representations of $\g$, the set $P_\ell$ represents the equivalence classes of representations of level at most $\ell$, meaning those representations $V_\lambda$ of $\g$ where $X^{\ell+1}$ acts trivially on $V_\lambda$ for every nilpotent element $X \in \g$.

\begin{rmk} \label{rmk-PellGamma} We note that the action of $\Gamma$ on $\g$ induces an action of $\Gamma$ on $P_\ell$ in the following way. Let $\phi \colon \g \times V \to V$ be a repreentation of $\g$, then we define the representation $\phi_{\gamma} \colon \g \times V \to V$ as $\phi_{\gamma}(X,v):=\phi(\gamma^{-1}X)v$ for all $X \in \g$ and $v \in V$. If the representation $\phi \in P_\ell$, then also $\phi_\gamma$ belongs to $P_\ell$ since $\Gamma$ sends nilpotent elements to nilpotent elements.\end{rmk}

The properties of $\tHh{V}$, for $V \in P_\ell$ are well known and described for example in \cite{kac1994infinite}, \cite{rainakac1988lecture}, \cite{tsuchiya1989conformal} and in \cite{beauville1994conformal}. The main results are collected in the following proposition. 

\begin{prop} \label{prop-HhnV} For $V \in P_\ell$ the following holds.
\begin{enumerate}
\item \label{uniquemaxquot} The module $\tHh{V}$ contains a maximal proper $U\HgL$ submodule $\Z_V$, so that it has a unique maximal irreducible quotient $\Hhn{V}:=\tHh{V}/\Z_V$. 
\item \label{VinHV} The natural map $V \to \Hhn{V}$ sending $v$ to $1 \otimes v$ identifies $V$ with the submodule of $\Hhn{V}$ annihilated by $UF^1\gL=U\g tk[\![t]\!]$.
\item \label{HhnVintegr} The module $\Hhn{V}$ is integrable, i.e. for any nilpotent element $X \in \g$ and every $f(t) \in k(\!(t)\!)$, the element $Xf(t)$ acts locally nilpotently on $\Hhn{V}$. \end{enumerate} \end{prop}

It follows that to every $(\Xt \to X \to \Spec(k), \tau) \in \bHurs{1}{\xi}(\Spec(k))$ and $V \in P_\ell$, we can associate the irreducible $U\HhL$ module $\Hhn{V}$ realized as the maximal irreducible quotient of $\tHh{V}$.

Let $(\Xt \to X \to S, \tau) \in \bHurs{1}{\xi}(S)$ and call $\sigma$ the composition $p\tau$. An isomorphism of $\hL$ with $\gL$ is fixed by $\tau$, as well as an isomorphism of ${\sigma}^*\h$ with $\g\otimes_k \Oo_S$. Denote by $\V:= V \otimes_k \Oo_S$ the extension of scalars of $V$ from $k$ to $\Oo_S$, so that $\V$ is naturally a representation of $\g \otimes_k \Oo_S={\sigma}^*\h$. We show how to construct $\Hhn{\V}$ as quotient of $\tHh{\V}$.

Let us assume first that $\Op \cong \Oo_S[\![t]\!]$, which provides an isomorphism $\tHh{\V} \cong \tHh{V} \otimes_k \Oo_S$ of $\HgL$ modules. Observe that this isomorphism does not depend on the choice of the parameter $t$ but only on the isomorphism $\hL \cong \gL$ 

It follows that $\tHh{\V}$ has a unique maximal $U \widehat{\g}\otimes_k k(\!(t)\!)$ proper submodule $\Z_V:=\Z_V \otimes \Oo_S$, where $\Z_V$ is the maximal proper submodule of $\Hhn{V}$. We define $\Hhn{\V}$ as the quotient $\tHh{\V} \Qdx \Z_S$ or equivalently as $\Hhn{V} \otimes \Oo_S$. This construction uses a choice of the isomorphism $\hL \cong \gL$, but since $\Z_V$ and hence $\Z_S$ satisfy a maximality condition, they do not depend on the isomorphism $\hL \cong \gL$, concluding that that $\Hhn{\V}$ is the maximal irreducible quotient representation of $U\HhL$ attached to $\V$.

We now drop the assumption $\Op \cong \Oo_S[\![t]\!]$. Since $\mathcal{I}_\sigma$ is locally principal, we can find an open covering $\{U_i\}$ of $X$ such that $\mathcal{I}_\sigma|_{U_i}$ is principal. This implies that $\varprojlim_n \Oo_{U_i}/ \left(\mathcal{I}_\sigma|_{U_i}\right)^n$ is isomorphic to $\Oo_{S_i}[\![t]\!]$ where $S_i:=\sigma^{-1}(U_i)$. Observe that this does not imply that $\Op \otimes_S \Oo_{S_i} \cong \Oo_{S_i}[\![t]\!]$, but only that $\Op \,\widehat{\otimes}_S \, \Oo_{S_i} \cong \Oo_{S_i}[\![t]\!]$. Consider then the sheaf of Lie algebras $\gL_i:=\g \otimes \Oo_{S_i}(\!(t)\!) \cong \gL \,\widehat{\otimes}\, \Oo_{S_i}$, and construct the $U\HgL_i$-module $\tHh{\V}_i$ as explained in the previous section.

\begin{claim} The inclusion $\gL \otimes \Oo_{S_i} \to \gL_i$ induces an isomorphism of $\Oo_{S_i}$-modules $\tHh{\V} \otimes_S \Oo_{S_i} \cong \tHh{\V}_i$.\end{claim}
\proof  We need to prove that $\tHh{\V} \otimes \Oo_{S_i} \cong \tHh{\V}_i$ is surjective. We use induction on the length of the elements of $U\gL_i$, where the length of an element $u \in U\gL_i$ is the minimum $n$ such that $u \in \oplus_{j=0}^n {\gL_i}^{\otimes j}$. Let $aX \in \gL_i$ with $X \in \g$ and $a=\sum_{i\geq -N} a_i t^i \in \Oo_{S_i}(\!(t)\!)$, and take $v \in \V$. The class of $aX \otimes v$ in $\tHh{\V}_i$ is the same as the one of $[aX]\otimes v :=[a]X \otimes v$, where $[a]=\sum_{i\geq -N}^{0} a_i t^i$, which then belongs to $\tHh{\V} \otimes_S \Oo_{S_i}$. Let now $Y=Y_1 \circ \dots \circ Y_n$ be an element of $U\gL_i$, and note that in $\tHh{\V}_i$ the element $Y \otimes v$ is equivalent to the class of $([Y_n]\circ \dots \circ [Y_1] + u) \otimes v$ where $u$ has length lower than $n$. Using the induction hypothesis we conclude the proof. \endproof

We define the $\Oo_{S_i}$-module $\Hhn{\V}|_{S_i}$ to be $\Hhn{\V}_i=\tHh{\V}_i \Qdx \Z_i$. This gives rise to the $\Oo_S$-module $\Hhn{\V}$ because on the intersection $S_{ij}$ the modules  $\Hhn{\V}_i$ and $\Hhn{\V}_j$ are isomorphic via to the transition morphisms defining $\mathcal{I}_\sigma$. Equivalently we could have defined $\Z|_{S_i}$ to be the image of $\Z_i$ in $\tHh{\V}|_{S_i}$ and so $\Hhn{\V}|_{S_i}$ would be the quotient of $\tHh{\V}|_{S_i}$ by $\Z|_{S_i}$. The modules $\Z|_{S_i}$ glue and give rise to a $\HgL$-module $\Z$ on $S$, so that $\Hhn{\V}$ is given by $\tHh{\V} \Qdx \Z$. This construction is invariant under the action of $\Gamma$, hence it defines $\HhL$ as a $U\HhL$-module.

\subsubsection{Integrable representations of level $\ell$ on $\bHur$} \label{subsubsec-IntegrReprLevEllHur}

We show here how to descend $\Hhn{\V}$ from $\bHurs{1}{\xi}$ to $\bHur$, so let consider $(\Xt \overset{q}\to X \overset{\pi}\to S, \sigma) \in \Hur(S)$. The first issue is that, unless we choose an isomorphism between $\sigma^*\h$ and $\g$, we are not able to provide a representation of $\sigma^*\h$ associated to a $V \in P_\ell$. In fact, one obstruction to this, as we noticed in Remark \ref{rmk-PellGamma}, is that $\Gamma$ does not in general act trivially on $P_\ell$. Moreover, an isomorphism between $\sigma^*\h$ and $\g$ exists only étale locally on $S$, so we cannot expect to associate to an element $V \in P_\ell$ a module on $\Hur$. The following set is what replaces $P_\ell$.

\begin{defi} A representation $\V$ of $\sigma^*\h$ is said to be of level at most $\ell$ if for every nilpotent element $X$ of $\sigma^*\h$, then $X^{\ell+1}$ acts trivially on $\V$. Equivalently this means that locally étale we can identify $\V$ with $V \otimes \Oo_S$ for a representation $V \in P_\ell$. Define $\Repl(\sigma^*\h)$ or by abuse of notation only $\Repl(\sigma)$ or $\Repl$ to be the set of isomorphism classes of irreducible and finite dimensional representations $\V$ of $\sigma^*\h$ of level at most $\ell$. \end{defi}

The main step towards the definition of the sheaf of conformal blocks attached to $\V \in \Repl$ is the following result.

\begin{prop} \label{prop-DescendtHh} Let $\V \in \Repl(\sigma^*\h)$. Then there exists a unique maximal proper $U\HhL$ submodule $\Z$ of $\tHh{\V}$.
\end{prop}

\proof We show that the maximal proper submodule of $\tHh{\V}$ on $\bHurs{1}{\xi}$ descends along $\textbf{Forg}^1_1$ to the maximal proper submodule of $\tHh{\V}$ on $\bHur$. Recall that since $\sigma(S)$ does not intersect the branch locus of $q$, we can find an étale covering $S'\to S$ such that the pullback of $(\Xt \to X \to S, \sigma)$ lies in the image of $\textbf{Forg}^1_1$. This implies that to give $\V \in \Repl$ is equivalent to give an irreducible and finite dimensional representation $\V'$ of ${\sigma'}^*\h$ and an isomorphism $\phi \colon p_1^*\V' \to p_2^*\V'$ satisfying the cocycle conditions on $S'''$, where $p_i \colon S''=S'\times_S S' \to S'$ is the projection on the $i$-th component. 

This tells us moreover that $\tHh{\V}$ is obtained by descending $\tHh{\V'}$ from $S'$ to $S$. Observe that up to the choice of an isomorphism ${\sigma'}^*\h \cong \g \otimes_k \Oo_{S'}$, the representation $\V'$ belongs to $P_\ell$, hence $\Z'$ and $\Hhn{\V'}$ are well defined. 

Since $\hL$ is a module on $\Oo_S$, we have a canonical isomorphism $\phi_{12} \colon {p_1}^*\hL|_{S'} \to {p_2}^*\hL|_{S'}$ satisfying the cocycle conditions on $S''':=S''\times_S S'$. Recall moreover that $\Z'$ is the maximal proper $U\HhL$ submodule of $\tHh{\V'}$, which is then $\Gamma$-invariant. This induces an isomorphism between $p_1^*\Z'$ and $p_2^*\Z'$ which satisfied the cocycle condition on $S'''$ and it is independent of the isomorphism $\hL \cong \gL$. This implies that $\Z'$ descends a $\HhL$-module $\Z$ on $S$ which is maximal by construction.
\endproof

For every $\V \in \Repl$, the maximal irreducible quotient of $\tHh{\V}$ by $\Z$ is denoted $\Hhn{\V}$ and, in view of Proposition \ref{prop-HhnV} satisfies the following properties:

\begin{cor} \label{cor-propHhn}
\begin{enumerate}
\item \label{VinHVtw} The natural map $\V \to \Hhn{\V}$ sending $v$ to $1 \otimes v$ identifies $\V$ with the submodule of $\Hhn{\V}$ annihilated by $UF^1\HhL$.
\item \label{HhnVintegrtw} The module $\Hhn{\V}$ is integrable.
\end{enumerate}
\end{cor}

Since $\hA$ is a Lie subalgebra of $\HhL$, it acts on $\Hhn{\V}$ by left multiplication.

\begin{defi} \label{def-HhnV} The \emph{sheaf of covacua} attached to $\V$  is the sheaf of $\Oo_S$ modules defined as \[\VV_\ell(\V)_{\Xt \to X}:=\hA \circ \Hhn{\V} \Qsx \Hhn{\V} = \hA \circ \Hhn{\V} \Qsx \tHh{\V} \Qdx \Z. \] By abuse of notation we are going to denote this $\Oo_S$-module simply by $\VV_\ell(\V)_{X}$. The \emph{sheaf of conformal blocks} attached to $\V$ is the dual $\VV_\ell(\V)_{X}^\dagger$ of the sheaf of covacua.

When $\V$ is the trivial representation of $\sigma^*\h$, we denote $\Hhn{\V}$ by $\Hhne$ and $\hA \Qsx \Hhne$ by $\VV_\ell(0)_{X}$.
\end{defi}

Given compatible families $\{\V(\sigma)\}_{\{\Xt \to X \to S,\sigma\}}$ defining an element $\V$ of $\Repl(\sigma_{\univ})$, the collection $\VV_\ell(\V(\sigma))_X$ defines $\VV_\ell(\V)$ on $\bHur$ which is called the \emph{universal sheaf of covacua}. Its dual module is the \emph{universal sheaf of conformal blocks}.

\subsubsection{Conformal blocks on $\bHur[n]$} \label{subsec-SemiLocalCase}

We extend the notion of covacua, and hence of conformal blocks, to the case in which more points of $X$, and by consequence more representations, are fixed. This will allow in a second time to express the factorization rules and, in view of the propagation of vacua, to drop the assumption that in the case $n=1$ we can only work with irreducible curves. This is explained in the classical contest in the last paragraphs of  \cite[Section 3]{looijenga2005conformal}. 

Let $(\Xt \overset{q}{\to} X \overset{\pi}{\to} S, \sigma_1, \dots, \sigma_n)$ be an $S=\Spec(R)$ point of $\bHur[n]$. For all $i \in \{1, \dots , n\}$ we denote by $S_i$ the divisor of $X$ defined by $\sigma_i$ and by $\mathcal{I}_i$ its ideal of definition. We denote by $X^*$ the open complement of $S_1 \cup \cdots \cup S_n$ in $X$ and we denote by $\hA$ the pushforward to $S$ of $\h$ restricted to $X^*$, in other words $\hA:= \pi_*(\h \otimes_{\Oo_X} \Oo_{X^*})$. As in the case $n=1$, we assume that $X^* \to S$ is affine.

In the same way as we defined $\Op$ in the case $n=1$, we set now $\Op_i$ to be the formal completion of $\Oo_X$ at $S_i$, i.e. $\Op_i=\pi_*\varprojlim_n \Oo_X / (\mathcal{I}_i)^n$. We set $\Ll_i= \varinjlim_{N} \pi_* \varprojlim_{n} \mathcal{I}_i^{-N} / \mathcal{I}_i^n$ and \[\hLi:=\varinjlim_N \pi_* \varprojlim_n \mathcal{I}_i^{-N} /\mathcal{I}_i^n \otimes_{\Oo_X} \h\] for all $i \in \{1, \dots, n\}$. The direct sum $\hLi[1] \oplus \cdots \oplus \hLi[n]$ is denoted by $\hL$ and $\Ll=\oplus \Ll_i$.

We extend centrally $\hLi$ in the same way as we did in the case $n=1$ obtaining $\HhLi$ with central element $c_i$. We denote by $\HhL$ the direct sum of $\HhLi$ modulo the relation that identifies all the central elements $c_i$'s so that \[ 0 \to c \Oo_S \to \HhL \to \hL \to 0\] is exact. The Lie algebra $\hA$ is still a sub Lie algebra of $\HhL$.

Let $i \in \{ 1, \dots, n\}$. We denote by $\Repl(i)$ the set of irreducible and finite dimensional representations of ${\sigma_i}^*\h$ of level at most $\ell$. As we have just done in Section \ref{subsubsec-IntegrReprLevEllHur} we attach to any $\V_i \in \Repl(i)$ the irreducible $U \HhLi$-module $\Hhn{\V_i}$. Taking their tensor product we obtain
\[\Hhn{\V_1, \dots, \V_n} := \Hhn{\V_i} \otimes \cdots \otimes \Hhn{\V_n}
\]
which then is a $U\HhL$-module with central charge $c$ acting by multiplication by $\ell$. Since $\hA$ is a Lie subalgebra of $\HhL$, we can take the sheaf of coinvariants with respect to that action. 

\begin{defi} The \emph{sheaf of covacua attached to $(\V_i)_{i=1}^n$} is the $\Oo_S$-module
\[\VV_\ell(\V_1, \dots, \V_n)_X := \hA \circ \Hhn{\V_1, \dots, \V_n} \Qsx \Hhn{\V_1, \dots, \V_n}.
\]The dual module $\VV_\ell(\V_1, \dots, \V_n)_X^\dagger$ is the \emph{sheaf of conformal blocks attached to $(\V_i)_{i=1}^n$}. \end{defi}

For every $i \in \{1, \dots, n\}$, we consider $\V_i$ as a representation of ${\sigma_{i,un}}^*\h$ defined by a compatible family $\{\V_i(\sigma_i)\}_{\{\Xt \to X \to S, \{\sigma_j\} \}}$ of representations of ${\sigma_i}^*\h$. The collection of $\VV_\ell(\V_1(\sigma_1), \dots, \V_n(\sigma_n))_X$ defines the module $\VV_\ell(\V_1, \dots, \V_n)$ on $\bHur[n]$ which we call the \emph{universal sheaf of covacua attached to $\{\V_i\}$}. Its dual module $\VV_\ell(\V_1, \dots, \V_n)^\dagger$ is the \emph{universal sheaf of conformal blocks attached to $\{\V_i\}$}

The main result of this paper is the following:

\begin{thm} \label{thm-LocFree} The sheaf of covacua $\VV_\ell(\V_1, \dots, \V_n)$  on $\bHur[n]$ is a vector bundle of finite rank. 
\end{thm}

As a first step in the direction of the proof, and inspired by \cite[Section 2.5]{sorger1996formule} we prove the following statement.

\begin{prop} \label{prop-HhnFinGen}  The sheaf $\VV_\ell(\V_1, \dots, \V_n)$ is a coherent module on $\bHur$. \end{prop}

\proof It is enough to show that the $\Oo_S$-module $\VV_\ell(\V_1, \dots, \V_n)_X$ is coherent and we show it in the case $n=1$. We furthermore observe that this is essentially a consequence of \cite[Lemma 2.5.2]{sorger1996formule}. As this is a local statement, we can assume that $\Ll \cong R(\!(t)\!)$ and we can fix an isomorphism $\hL \cong \gL$. Observe that the quotient $\hA \Qsx \hL \Qdx F^0\hL$ is a finitely generated $R$-module as it computes $H^1(X,\h)$ and $\h$ is locally free over $X$. This implies that $\hA \Qsx \HhL \Qdx F^1\hL$ is finitely generated too over $R$ and so we can choose finitely many generators $e_1, \dots, e_n$ so that we can write \[ \HhL=F^1\hL + \hA + \sum_{i=1}^n R e_i.\]
which in terms of enveloping algebras becomes
\[U\HhL = \sum_{(N_1, \dots, N_n) \in {\ZZ_{\geq 0}}^n} U(\hA) \circ {e_1}^{\circ N_1} \circ \cdots \circ {e_n}^{\circ N_n} \circ U(F^1\hL)
\] thing that can be proven using induction on the length of elements of $U\HhL$.

We can furthermore assume that the elements $e_i$ acts locally nilpotently on $\Hhn{\V}$, meaning that there exists $M \in \ZZ_{> 0}$ such that $e_i^{\circ M}$ acts trivially on $\Hhn{\V}$. In fact we might use the isomorphism $\hL$ with $\gL$ and the Cartan decomposition of $\g=\mathfrak{t} \oplus_{\alpha \in R(\g,\mathfrak{t})} \g_\alpha$. The algebras $\g_{\alpha}$'s are nilpotent and generate $\g$, so that $\gL$ is generated by $\oplus_{\alpha \in R(\g,\mathfrak{t})} \g_\alpha\Ll$. This means that also the elements $e_i$ are generated by elements of $\oplus_{\alpha \in R(\g,\mathfrak{t})}\g^\alpha\Ll$ so that, up to replace $e_i$ with a choice of nilpotent generators, we can ensure that all the $e_i$'s live in $\oplus_{\alpha \in R(\g,\mathfrak{t})}\g^\alpha\Ll$ and so using Corollary \ref{cor-propHhn} (\ref{HhnVintegrtw}) the $e_i$'s will act locally nilpotently on $\Hhn{\V}$.

It follows that \[\tHh{\V}= \sum_{(N_1, \dots, N_n) \in {\ZZ_{\geq 0}}^n} U(\hA) \circ {e_1}^{\circ N_1} \circ \cdots \circ {e_n}^{\circ N_n} \underset{c=\ell}\otimes \V  
\]
and that
\[\Hhn{\V} = \sum_{(N_1, \dots, N_n) \in {\ZZ_{\geq 0}}^n} U(\hA) \circ {e_1}^{\circ N_1} \circ \cdots \circ {e_n}^{\circ N_n} \underset{c=\ell}\otimes \V \Qdx \Z
\]
Using induction on $n$ and the fact that the $e_i$'s act locally nilpotently, we can conclude that the sum can be taken over finitely many $(N_1,\dots,N_n) \in {\ZZ_{\geq 0}}^n$, hence that the quotient $\hA \Qsx \Hhn{\V} = \VV_\ell(\V)_X$ is finitely generated.\endproof

\section{The projective connection on $\VV_\ell(\V)$} \label{sec-ProjConn}

We want to prove that the sheaf of covacua $\VV_\ell(\V_1, \dots, \V_n)$ is a vector bundle on the Hurwitz stack $\bHur[n]$. This will in particular imply that also its dual $\VV_\ell(\V_1, \dots, \V_n)^\dagger$ is a vector bundle, and that its rank is constant on $\bHur[n]$. Since we already know that $\VV_\ell(\V_1, \dots, \V_n)$ is coherent, one method to exhibit local freeness is to provide a projectively flat connection on it. In this section we provide a projective action of $\T_{\bHur / k}(-\log(\Delta))$ on $\VV_\ell(\V)$, showing its freeness when restricted to $\Hur$. We will explain in detail how to achieve this in the case $n=1$, and postpone to the end of the section the situation with more marked points.

\subsection{The tangent to $\bHur$}  \label{subsec-tangHur} Let $(\Xt \to X, \sigma) \in \bHur(\Spec(k))$ and recall that in Remark \ref{rmk-NCD} we saw that the tangent space of $\bHur$ at $(\Xt \to X, \sigma)$ is isomorphic to the tangent space of $\Mgn[(1+\deg(\xi))]$ at $(X, \sigma \cup \Rr)$. 
The latter, which is the space of infinitesimal deformations of $(X, \sigma \cup \Rr)$, can be explicitely described as the space $\Ext^1(\Omega_{X/k},  \Oo(-\Rr-\sigma(S))$ \cite[Chapter XI]{ArbarelloCornalba} which sits in the short exact sequence
\[ 0 \to H^1(X, \T_{X/k}(-\Rr -\sigma(S)) \to \Ext^1(\Omega_{X/k},  \Oo(-\Rr-\sigma(S)) \to H^0( X, \Extt^1(\Omega_{X/k},  \Oo(-\Rr-\sigma(S))) \to 0
\] where the last term is supported on the singular points of $X$.

We now use the assumption that $(X, \sigma \cup \Rr)$ is a stable marked curve to ensure that there exists a versal family $\Xx \to S$ with a reduced divisor $\sigma_\Xx + \Rr_{\Xx}$ deforming it and such that the subscheme of $S$ whose fibres are singular is a normal crossing divisor $\Delta$. Call $s_0$ the point of $S$ such that $\Xx|_{s_0}$ is $X$. The versality condition means that the Kodaira-Spencer map
\[ \texttt{KS} \colon \T_{S/k} \to \Extt^1(\Omega_{\Xx/S}, \Oo(-\Rr_{\Xx}-\sigma_{\Xx}))\] is an isomorphism, so that we identify the tangent space of $\bHur$ at $(\Xt \to X, \sigma)$ with the tangent space of $S$ at $s_0$. The conclusion is that to provide a projective connection on $\VV_\ell(\V)$ is equivalent to provide an action of $\T_{S/k}$ on $\VV_\ell(\V)_{X}$ for every versal family $\Xt \overset{q}\to X \overset{\pi}\to S$. As aforementioned, we will however not be able to provide a projective action of the whole $\T_{S/k}$, but only of the submodule $\T_{S/k}(-\log(\Delta))$, which via the Kodaira-Spencer map is identified with $R^1\pi_*(\T_{X/S}(-\Rr-\sigma(S)))$.

\subsection{Tangent sheaves and the action of $\Gamma$} 	\label{subsec-tangent}

In view of the previous observations, we assume that $(\Xt \overset{q}\to X \overset{\pi}\to S, \sigma) \in \bHur(S)$ is a versal family, so that the locus of points $s$ of $S$ such that the fibres $X_s$ (or equivalently $\Xt_s$) are non smooth is a normal crossing divisor $\Delta$ of $S$. We give in this section a description of $\T_{S/k}(-\log(\Delta))$ and $\T_{S/k}$, by realizing it as a quotient of certain sheaves of derivations. We already introduced in Subsection \ref{subsubsec-ResiduePairing} the $\Oo_S$-module $\theta_{\Ll/S}$ of continuous $\Oo_S$-linear derivations of $\Ll$ and we define now $\theta_{\Ll,S}$ as the $\Oo_S$-module of continuous $k$-linear derivations of $\Ll$ which restrict to derivations of $\Oo_S$. Observe that $\theta_{\Ll/S}$ and $\theta_{\Ll,S}$ depend only on the marked curve $(X \to S, \sigma)$, so the following well known result belongs to the classical setting.

\begin{prop} \label{prop-exseq} The sequence of $\Oo_S$-modules 
\[ 0 \to \theta_{\Ll/S} \to \theta_{\Ll,S} \to \T_{S/k} \to 0
\] is exact.
\end{prop}

In similar fashion we now describe the subsheaf $\T_{S/k}(-\log(\Delta))$ as quotient of appropriate sheaves of derivations. Following the notation of \cite{looijenga2005conformal}, we denote by $\theta_{\Aa/S}$ the sheaf of derivations $f_*\T_{\Xt^*/S}$ and in view of Proposition \ref{prop-tang}, we write $\theta_{\Aa/S}(-\Rr)$ to denote $(f_*(\T_{\Xt^*/S}))^\Gamma$. In a similar way we consider the action of $\Gamma$ on the pushforward to $S$ of $\T_{\Xt^*,S}$, the sheaf of $k$-linear derivations of $\Oo_{\Xt^*}$ which restrict to derivations of $f^{-1}\Oo_S$ and we call $\theta_{\Aa,S}(-\Rr)$ the sub module of $\Gamma$-invariants.

\begin{rmk} Recall that we defined $\widetilde{\Ll}$ as $\varinjlim_{N \in \ZZ_{\geq 0}} f_* \varprojlim_{n \in \ZZ_{> 0}} (q^*\mathcal{I}_\sigma)^{-N} / (q^*\mathcal{I}_\sigma)^n$ and define now \[ \theta_{\widetilde{\Ll}/S} := \varinjlim_{N \in \ZZ_{\geq 0}} \pi_*q_* \varprojlim_{n \in \ZZ_{> 0}} \T_{\Xt/S} \otimes_{\Oo_\Xt} q^*\mathcal{I}_\sigma^{-N} / q^*\mathcal{I}_\sigma^n\]
or equivalently $\theta_{\widetilde{\Ll}/S}$ is the $\Oo_S$-module of continuous derivations of $\widetilde{\Ll}$ which are $\Oo_S$-linear. Thanks to Proposition \ref{prop-tang}, the $\Oo_S$-submodule of $\Gamma$-invariants of $\theta_{\widetilde{\Ll}/S}$ is identified with $\varinjlim_{N \in \ZZ_{\geq 0}} \pi_* \varprojlim_{n \in \ZZ_{> 0}} \T_{X/S}(-\Rr) \otimes_{\Oo_X} \mathcal{I}_\sigma^{-N} / \mathcal{I}_\sigma^n$ which equals $\varinjlim_{N \in \ZZ_{\geq 0}} \pi_* \varprojlim_{n \in \ZZ_{> 0}} \T_{X/S}\otimes_{\Oo_X} \mathcal{I}_\sigma^{-N} / \mathcal{I}_\sigma^n$ as $\Rr$ and $\sigma(S)$ are disjoint. The latter is the $\Oo_S$-module of continuous and $\Oo_S$-linear derivations of $\Ll$, which is $\theta_{\Ll/S}$.
\end{rmk} 

The previous remark implies moreover that $\theta_{\Aa/S}$ is a submodule of $\theta_{\Ll/S}$. 

\begin{rmk} \label{rmk-coeffder} Observe that the action of $\T_{\Xt/S}$ (resp. of $\T_{\Xt,S}$) on $\g \otimes_k \Oo_\Xt$ by coefficientwise derivation is $\Gamma$-equivariant. This implies that $(q_*\T_{\Xt/S})^\Gamma$ (resp. $(q_*\T_{\Xt,S})^\Gamma$) acts on $\h$ and we will say that the action is by coefficientwise derivation. In particular $\theta_{\Aa/S}(-\Rr)$ (resp. $\theta_{\Aa,S}(-\Rr)$) acts on $\hA$ by coefficientwise derivation and the same holds for $\theta_{\Ll/S}$ and $\theta_{\Ll,S}$ acting on $\hL$.
\end{rmk}

\begin{prop} \label{prop-exseqALog} The sequence \[ 0 \to \theta_{\Aa/S}(-\Rr) \to \theta_{\Aa,S}(-\Rr) \to \T_{S}(-\log(\Delta)) \to 0
\]is exact.
\end{prop}

\proof As taking $\Gamma$-invariants is an exact functor ($\text{char}(k)=0$) and $\Gamma$ acts trivially on $\T_{S/k}(-\log(\Delta))$, it suffices to prove that the sequence 
\[0 \to f_*\T_{\Xt^*/S} \to f_*\T_{\Xt^*,S} \to \T_S(-\log(\Delta)) \to 0 \]is exact. This statement does not depend on the covering, and appears in \cite{sorger1996formule} and \cite{looijenga2005conformal}.
\endproof

\subsection{The Virasoro algebra of $\Ll$} \label{subsec-VirasoroAlgebra}

Now that we can express $\T_{S/k}(-\log(\Delta))$ as $\theta_{\Aa,S}(-\Rr) \Qdx \theta_{\Aa/S}(-\Rr)$, we will define a projective action of $\theta_{\Aa,S}(-\Rr)$ on $\VV_\ell(\V)$ which factors through that quotient.

In order to achieve this result, we will follow the methods of \cite{looijenga2005conformal} and define as first step the Virasoro algebra $\widehat{\theta}_{\Ll/S}$ of $\Ll$ as a central extension of $\theta_{\Ll/S}$. We report in this section the construction of $\widehat{\theta}_{\Ll/S}$ as explained in \cite[Section 2]{looijenga2005conformal}, for which we will use the same notation. As mentioned before, the $\Oo_S$-module $\theta_{\Ll/S}$ does not depend on the covering, and the same holds for its central extension $\widehat{\theta}_{\Ll/S}$. The reader who is already familiar with the construction, can therefore skip this section.  

\subsubsection*{The Lie algebras $\lL$ and its central extension $\widehat{\lL}$}

We denote by $\lL$ the sheaf of abelian Lie algebras (over $S$) whose underlying module is $\Ll$. The filtration $F^*\!\Ll$ gives the filtration $F^*\lL$. Denote by $U\lL$ the universal enveloping algebra, which is isomorphic to $\Sym(\lL)$ since $\lL$ is abelian. This algebra is not complete with respect to the filtration $F^*\lL$, so we complete it on the right obtaining \[\overline{U}\lL:= \varprojlim_n U\lL / U\lL \circ F^n\lL.\]

\begin{rmk} Note that in this case the completions on the right $\varprojlim_n U\lL / U\lL \circ F^n\lL$ and on the left $\varprojlim_n F^n\lL \circ U\lL \Qsx U\lL$ coincide because $\lL$ is abelian.  The element $\sum_{i\in\ZZ_{> 0}} t^{-i} \circ t^i$ belongs to $\overline{U}\lL$, as well as $\sum_{i\in\ZZ_{> 0}} t^{-(i)^m} \circ t^i$ for every $m \in \ZZ_{> 0}$. However $\sum_{i \in \ZZ_{> 0}}t^{-i} \circ t$ is not an element of $\overline{U}\lL$. \end{rmk}

We extend centrally $\lL$ via the residue pairing described in Section \ref{subsubsec-ResiduePairing} defining the Lie bracket on $\widehat{\lL}=\lL \oplus \hslash\Oo_S$ as \[ [f + \hslash r, g + \hslash s] = \hslash \Res(gdf)\]
for every $f,g \in \lL$ and $r,s \in \Oo_S$. The filtration of $\lL$ extends to a filtration of $\widehat{\lL}$ by setting $F^i\widehat{\lL}=F^i\lL$ for $i \geq 0$ and $F^i\widehat{\lL}=F^i\lL \oplus \hslash \Oo_S$ for $i \leq 0$. The universal enveloping algebra of $\widehat{\lL}$ is denoted by $U\widehat{\lL}$ and $\overline{U} \widehat{\lL}$ denotes its completion on the right with respect to the filtration $F^*\widehat{\lL}$. Note that since $\hslash$ is a central element, we have that $\overline{U} \widehat{\lL}$ is an $\Oo_S[\hslash]$ algebra so that we will write $\hslash^2$ instead of $\hslash \circ \hslash$ and similarly $\hslash^n$ for every $n \in \ZZ_{> 0}$.

\begin{rmk} Since $\widehat{\lL}$ is no longer abelian, completion on the right and on the left differ. Take for example the element $\sum_{i \in \ZZ_{> 0}} t^i \circ t^{-i}$ which belongs to $\overline{U}\lL$. It does not belong to $\overline{U\widehat{\lL}}$: an element on the completion on the right morally should have zeros of increasing order on the right side, but in this case, in order to "bring the element $t^i$ on the right side", we should use the equality $t^i \circ t^{-i}=t^{-i} \circ t^i + \hslash i$ infinitely many times, which is not allowed.\end{rmk}

\subsubsection*{The Virasoro algebra of $\Ll$} \label{subsubsec-Virasoro}

We want to use the residue morphism $\Res \colon \omega_{\Ll/S} \to \Oo_S$ to view $\theta_{\Ll/S}$ as an $\Oo_S$ submodule of $\overline{U}\widehat{\lL}$ and induce from this a central extension. Let $D \in \theta_{\Ll/S}$, and since $\omega_{\Ll/S}$ and $\theta_{\Ll/S}$ are $\Ll$ dual we identify $D$ with the map
\[ \phi_D \colon \omega_{\Ll/S} \times \omega_{\Ll/S} \to \Oo_S, \quad (\alpha, \beta) \mapsto \Res(D(\alpha)\beta).
\]
Notice that since $\Res(D(\alpha)\beta)=\Res(D(\beta)\alpha)$, we have that $\phi_D$ is an element of $(\Sym^2(\omega_{\Ll/S}))^\vee$. Moreover, since the latter is canonically isomorphic to the closure of $\Sym^2(\lL)$ in $\overline{U}\lL$, we will consider $\phi_D$ as an element of $\overline{U}\lL$. We define $C \colon \theta_{\Ll/S} \to \overline{U}\lL$ by setting $2C(D)=\phi_D$.

\begin{rmk} \label{rmk-explicitCDF} Assume for simplicity that $R=k$ and identify $\Ll$ with $k(\!(t)\!)$. For every $i \in \ZZ$ we set $\alpha_i=t^{-i-1}dt$ and $a_i=t^{i}$ so that $\Res(a_i, \alpha_j)=\delta_{ij}$ and $\{\alpha_i\}$ and $\{a_i\}$ are linearly independent generators of $\omega_{\Ll/S}$ and $\Ll$. Then we can write explicitly
\[C(D)=\dfrac{1}{2}\sum_{i\in \ZZ} D(t^{-i-1}dt) \circ t^i .
\]
In general, let $\{\alpha_i\}$ and $\{a_i\}$ be linearly independent generators of $\omega_{\Ll/S}$ and $\Ll$ with the property that $\Res(a_i, \alpha_j)=\delta_{ij}$. Then we can write \[C(D)=\dfrac{1}{2}\sum_{i \in \ZZ}  D(\alpha_i) \circ a_i\] which is a well defined object of $\overline{U}\lL$ thanks to the previous remark. 
\end{rmk}

As explained in \cite[Section 2]{looijenga2005conformal}, the central extension $\widehat{\lL}$ and the inclusion $C \colon \theta_{\Ll/S} \to \overline{U}\widehat{\lL}$ induce a central extension $\widehat{\theta}_{\Ll/S}$ of $\theta_{\Ll/S}$. We recall here how this is achieved. Consider now $\lL \otimes \lL$, and call $\lL_2$ its image in $U\widehat{\lL}$. This means that $\lL_2=\lL \otimes \lL \oplus \hslash \Oo_S$ modulo the relation $f \otimes g= g\otimes f + \hslash \Res(gdf)$. Denote by $\overline{\lL_2}$ its closure in $\overline{U}\widehat{\lL}$ and observe the following diagram
\[\xymatrix{ 0 \ar[r] &
\hslash \Oo_S \ar[r] & \overline{\lL_2} \ar[r]^-{[-]_\hslash} & \overline{\Sym^2(\lL)} \ar[r] & 0\\
&&& \theta_{\Ll/S} \ar[u]_C}
\]
where $[-]_\hslash$ is the reduction modulo the central element $\hslash \Oo_S$ so that the short sequence is exact.

\begin{defi} We define $\widehat{\theta}_{\Ll/S}$ to be the pullback of $\theta_{\Ll/S}$ along $[-]_\hslash$. Equivalently its elements are pairs $(D,u) \in \theta_{\Ll/S} \times \overline{\lL_2}$ such that $C(D) = u \mod \hslash \Oo_S$. 
\end{defi}

Denote by $\widehat{C} \colon \widehat{\theta}_{\Ll/S} \to \overline{U}\widehat{\lL}$ the injection $\widehat{C} (D,u)=u$ and we write $[-]^\theta_\hslash$ for the pullback of $[-]_\hslash$ along $C$ so that we have the commutative diagram with exact rows
\[\xymatrix{ 0 \ar[r] &
\hslash \Oo_S \ar[r] & \overline{\lL_2} \ar[r]^-{[-]_\hslash} & \overline{\Sym^2(\lL)} \ar[r] & 0\\
0 \ar[r] & \hslash \ar@{=}[u]\Oo_S \ar[r] & \widehat{\theta}_{\Ll/S} \ar[u]_{\widehat{C}} \ar[r]\ar[r]^{[-]^\theta_\hslash} & \theta_{\Ll/S} \ar[u]_C \ar[r] &0.} 
\]

Observe, for example using Remark \ref{rmk-explicitCDF}, that the map $C$ is not a Lie algebra morphism, and so $\widehat{\theta}_{\Ll/S}$ does not arise naturally as a Lie algebra which centrally extends $\theta_{\Ll/S}$. 

We however want to induce a Lie bracket on $\widehat{\theta}_{\Ll/S}$ from the one of $\overline{U}\widehat{\lL}$ by conveniently modifying $\widehat{C}$. To understand how to do this, local computations are carried out.

\begin{defi} \label{def-ord} Choose a local parameter $t$ so that locally $\Ll \cong \Oo_S(\!(t)\!)$. Define the \emph{normal ordering} $: \star : \colon \overline{\Sym^2(\lL)} \to \overline{\lL^2}$ by setting 
\[: t^n \otimes t^m: = \begin{cases} t^n \otimes t^m &n \leq m\\ t^m \otimes t^n &n \geq m
\end{cases}\] and extend it by linearity to every element of $\overline{\Sym^2(\lL)}$. \end{defi}

The map $: \star :$ defines a section of $[-]_\hslash$, so that $(\text{Id}, : \star : C)$ is a section of $[-]_\hslash^\theta$. Once we make the choice of a local parameter defining the ordering $:\star:$, we will denote by $\widehat{D}$ the element $(D, :C(D):) \in \widehat{\theta}_{\Ll/S}$. Consider the following relations which hold in $\overline{U}\widehat{\lL}$ and which are proved in \cite[Lemma 2.1]{looijenga2005conformal}.

\begin{lem} \label{lem-LemLie} Let $D \in \theta_{\Ll/S}$ and $D_i=t^{i+1}d/dt \in \theta_{\Ll/S}$, Then we have \begin{enumerate}
\item $[\widehat{C}(\widehat{D}), f]=-\hslash D(f)$ for every $f \in \lL \subset \widehat{\lL}$;
\item $[\widehat{C}(\widehat{D}_k), \widehat{C}(\widehat{D}_l)]=-\hslash (l-k) \widehat{C}(\widehat{D}_{k+l}) + \dfrac{k^3-k}{12} \hslash^2  \delta_{k,-l}$. \qed
\end{enumerate}
\end{lem}

This suggests to rescale the morphism $\widehat{C}$ and to define
\[ T:= -\dfrac{\widehat{C}}{\hslash} \colon \widehat{\theta}_{\Ll/S} \to \overline{U}\widehat{\lL}\left[ \dfrac{1}{\hslash}\right]
\]
which is injective and its image is a Lie subalgebra of the target. Denote by $c_0$ the element $(0, -\hslash)$ which is sent to $1$ by $T$. By construction we obtain the following result.

\begin{prop} \cite[Corollary-Definition 2.2]{looijenga2005conformal} The Lie algebra structure induced on $\widehat{\theta}_{\Ll/S}$ by $T$ is a central extension of the canonical Lie algebra structure on $\theta_{\Ll/S}$ by $c_0 \Oo_S$. This is called the \emph{Virasoro algebra} of $\Ll$.\end{prop}

\subsection{Sugawara construction} \label{subsec-Sugawara}

In this section we generalize to our case, i.e. using $\HhL$ in place of $\HgL$, the construction of \[T_\g \colon \widehat{\theta}_{\Ll/S} \to \left(\overline{U}\HgL[(c+\check{h})^{-1}]\right)^{\Aut(\g)}\] described by Looijenga in \cite[Corollary 3.2]{looijenga2005conformal}, which essentially represents the local picture of our situation. In the classical case the idea is to use the Casimir element of $\g$ to induce, from $\widehat{C}$ the map $\widehat{C}_\g \colon \widehat{\theta}_{\Ll/S} \to \overline{U}\HgL$ which, in turn, will give the map of Lie algebras $T_\g$. When in place of $\gL$ we have $\hL$, we can run the same argument using the element Casimir $\mathfrak{c}$ of $\hL$. 

As in Section \ref{subsec-CentralExtensionhL}, we consider the normalized Killing form defined on $\g$. Recall that it provides an isomorphism between $\g$ and $\g^\vee$, hence it gives an identification of $\g \otimes \g$ with $\text{End}_k(\g)=\g \otimes \g^\vee$. Moreover, as $\sigma(S)$ is disjoint from the ramification locus, we also have that $(\,|\,)_\hL$ provides an isomorphism of $\hL$ with $\hL^\vee$, giving in this way an identification of $\hL \otimes_\Ll \hL$ with $\text{End}_\Ll(\hL)$. The \emph{Casimir element} of $\hL$ with respect to the form $(\,|\,)_\hL$ is the element in $\hL \otimes_\Ll \hL$ corresponding to the identity $\text{Id}_\hL$ via the identification provided by $(\,|\,)$. We denote it by $\mathfrak{c}$.

\begin{rmk} We could have defined the Casimir element of $\hL$ via the local isomorphism of $\hL$ with $\gL$. Let $\mathfrak{c}(\g)$ be the Casimir element of $\g$, and observe that via the inclusion $\g \to \gL$, we can see it as an element of $\gL \otimes_\Ll \gL$. Since $\mathfrak{c}(\g)$ is invariant under automorphisms, it is invariant under $\Gamma$, hence it gives an element $\hL \otimes_\Ll \hL$ which equals $\mathfrak{c}$. \end{rmk}

Since $(\,|\,)$ is a symmetric form, we have that also $\mathfrak{c}$ is a symmetric element of $\hL \otimes_\Ll \hL$ and moreover $\mathfrak{c}$ lies in the centre of $U_\Ll(\hL)$. As $\hL$ is simple over $\Ll$, this implies that there exists $\check{h} \in k$ such that $\text{ad}(\mathfrak{c})X= 2\check{\h} X$ for all $X \in \hL$, where $\text{ad}(-)$ denotes the adjoint representation of $\hL$. 

\begin{rmk} \label{rmk-basCas} Locally, for every bases $\{X_i\}_{i=1}^{\dim(\g)}$ and $\{Y_i\}_{i=1}^{\dim(\g)}$ of $\hL$ such that $(X_i|Y_j)_\hL=\delta_{ij}$ we have the explicit description of $\mathfrak{c}$ as $\sum_{i=1}^{\dim(\g)} X_i \circ Y_i$. It follows that $\check{h}$ is given by the equality $\sum_{i=1}^{\dim(\g)}[ X_i, [Y_i, Z]]=2 \check{h} Z$ for every $Z \in \hL$.
\end{rmk}

Let denote by $\overline{U}\HhL$ the completion on the right of $U\HhL$ with respect to the filtration $F^*\HhL$ given by $F^n\hL$ for $n \geq 1$. We now construct $\widehat{\gamma}_{\mathfrak{c}} \colon \overline{\lL} \to \overline{U}\HhL$ which composed with $\widehat{C}$ will give $\widehat{C}_{\h} \colon \widehat{\theta}_{\Ll/S} \to \overline{U} \HhL$.

Let consider the map $\gamma_{\mathfrak{c}} \colon \lL \otimes_{\Oo_S} \lL \to \hL \otimes_{\Oo_S} \hL \subset U\HhL$ given by tensoring with $\mathfrak{c}$. This map uniquely extends to a map of Lie algebras $\lL_2 \to U\HhL$ as follows. Using local bases as in Remark \ref{rmk-basCas} and the symmetry of ${\mathfrak{c}}$ we deduce the following equality \[\gamma_{\mathfrak{c}}\left(\hslash\Res(gdf)\right)=\gamma_{\mathfrak{c}}(f \circ g - g \circ f)= c \dim(\g)\Res(gdf)+ c \sum_{i=1}^{\dim(\g)} \Res(fg(dY_i|X_i))\] and recalling that Remark \ref{rmk-dkill} implies that $(dX_i|Y_i)=-(dY_i|X_i)$, we conclude that
\[\gamma_{\mathfrak{c}}\left(\hslash\Res(gdf)\right)=c \dim(\g)\Res(gdf).\]
We then define $\widehat{\gamma}_{\mathfrak{c}} \colon \lL_2 \to U\HhL$ by sending $\hslash$ to $c \dim(\g)$ and acting as $\gamma_{\mathfrak{c}}$ on $\lL \otimes \lL$.  Such a map can be extended to the closure of $\lL_2$ in $\overline{U}\widehat{\lL}$ once we extend the target to $\overline{U}\HhL$, obtaining $\widehat{\gamma}_{\mathfrak{c}} \colon \overline{\lL_2} \to U\HhL$.
We define $\widehat{C}_\h$ as the composition \[\widehat{\gamma}_{\mathfrak{c}}  \widehat{C} \, \colon \widehat{\theta}_{\Ll/S} \to \overline{U} \HhL
.\]

As for $\widehat{C}$, also in this case the morphism $\widehat{C}_\h$ does not preserve the Lie bracket, and thanks to local computation we understand how to solve this issue. Following \cite{tsuchiya1989conformal} we first of all extend the normal ordering defined in \ref{def-ord} as as follows. 

\begin{defi} \label{def-circord} Let fix an isomorphism between $\Ll$ and $R(\!(t)\!)$ for a local parameter $t \in \mathcal{I}_\sigma$. Let $Xt^n$ and $Yt^m$ be elements of $\gL=\g \otimes R(\!(t)\!)$. Then we set
\[\circdot \,X t^n \otimes Y t^m \,\circdot = \begin{cases} Xt^n \otimes Yt^m &n < m\\ \dfrac{1}{2} 	(X t^n \otimes Y t^m + Y t^m \otimes X t^n)	&n=m	\\ Yt^m \otimes Xt^n &n > m.
\end{cases}\] This definition is $\Gamma$-equivariant, hence defines a normal ordering on $\hL \otimes \hL$. \end{defi} This defines a section from the image of $\gamma_c$ to the image of $\widehat{\gamma}_c$ which makes the diagram to commute:

\[\xymatrix{ 0 \ar[r] &
c \Oo_S \ar[r] & \text{Im}(\widehat{\gamma}_c) \ar[r]^-{[-]_c} & \text{Im} (\gamma_c) \ar@/^1.5pc/[l]^{\circdot\star\circdot}\ar[r] & 0\\
0 \ar[r] &
\hslash \Oo_S \ar[r] \ar[u] & \overline{\lL_2} \ar[r]^-{[-]_\hslash}\ar[u]_{\widehat{\gamma}_c} & \overline{\Sym^2(\lL)} \ar[r] \ar[u]_{\gamma_c} \ar@/^1.5pc/[l]^{:\star:}& 0\\
&& \widehat{\theta}_{\Ll/S} \ar[u]_{\widehat{C}} \ar[r]& \theta_{\Ll/S} \ar[u]_C}
\]

For any $D \in \theta_{\Ll/S}$, we write $\widehat{C}_\h(\widehat{D})$ to denote the element
$\circdot \, \gamma_\mathfrak{c} C (D)\, \circdot \, = \widehat{\gamma}_\mathfrak{c} : C(D) :$.

\begin{rmk} \label{rmk-explLoc} As we have done in Remark \ref{rmk-explicitCDF} we can write locally the element $\widehat{C}_\h(\widehat{D})$ in a more explicit way. Consider the morphism $1 \otimes (\,|\,) \colon \omega_{\Ll/S} \otimes \hL \to \omega_{\Ll/S} \otimes \hL^\vee$ and, after tensoring it with $\hL$, compose it with $\Res_\h$ to obtain the pairing $\Res_{(\,|\,)} \colon \omega_{\Ll/S} \otimes_\Ll \hL \times \hL \to \Oo_S$. Let $\{A_i\}$ and $\{B_i\}$ be orthonormal bases of $\omega_{\Ll/S} \otimes \hL$ and $\hL$ with respect to $\Res_{(\,|\,)}$. Then for every $D \in \theta_{\Ll/S}$ we have
\[ \widehat{C}_\h(\widehat{D})= \frac{1}{2} \sum \circdot D(A_i) \circ B_i \circdot
\]
where we see $D$ as a linear map $\omega_{\Ll/S} \to \Ll$, so that $D(A_i) \in \hL$.
\end{rmk}

As in \cite[Lemma 3.1]{looijenga2005conformal} we have the following result.

\begin{lem}\label{lem-Lemdue} The following equalities hold true in $\overline{U}\HhL$:
\begin{enumerate}
\item $[\widehat{C}_\h(\widehat{D}), X]=-(c+\check{h})D(X)$ for all $X \in \hL$ and $D \in \theta_{\Ll/S}$;
\item $[\widehat{C}_\h(\widehat{D}_k),\widehat{C}_\h(\widehat{D}_l)]=(c+ \check{h})(k-l) \widehat{C}_\h(\widehat{D}_{k+l}) + c \dim(\g)(c+\check{h}) \dfrac{k^3-k}{12}\delta_{k,-l}$ where $D_i=t^{i+1}d/dt$.
\end{enumerate}
\end{lem}

As Lemma \ref{lem-LemLie} also Lemma \ref{lem-Lemdue} suggests to rescale $\widehat{C}_\h$ and consider instead the map
\[T_\h := -\dfrac{\widehat{C}_\h}{c+\check{h}} \colon \widehat{\theta}_{\Ll/S} \to \overline{U}\HhL\left[ \dfrac{1}{c+\check{h}}\right]
\]
which is compatible with the Lie brackets of $\widehat{\theta}_{\Ll/S}$ and $\overline{U}\HhL[(c+\check{h})^{-1}]$, proving the following statement.

\begin{prop} \label{prop-Sugawara} The map $T_\h$ is a homomorphism of Lie algebras which sends the central element $c_0=(0, -\hslash)$ to $(c \dim(\g))/(c+ \check{h})$. We call $T_\h$ the \emph{Sugawara representation} of $\widehat{\theta}_{\Ll/S}$.
\end{prop}

\subsubsection{Fock representation}\label{subsubsec-fockg}

We induce the representation $T_\h$ to the quotient $\Ff^+(\hL)$ of $\overline{U}\HhL$ defined as \[\Ff^+(\hL):=\left( U\HhL \Qdx U\HhL \circ F^1\HhL \right) \left[ \dfrac{1}{c +\check{h}} \right] = \left( \overline{U}\HhL \Qdx \overline{U} \HhL \circ F^1\HhL \right) \left[ \dfrac{1}{c +\check{h}} \right]
\]
By abuse of notation call $T_\h$ the composition of $T_\h$ with the projection of $\overline{U}\HhL$ to $\Ff^+(\hL)$, so that $\Ff^+(\hL)$ is a representation of $\widehat{\theta}_{\Ll/S}$. 
We can depict the result as follows
\[\xymatrix{0 \ar[r] &\Oo_S \cdot \text{Id}  \ar[r]  &\End(\Ff^+(\hL)) \ar[r] & \PP\End(\Ff^+(\hL)) \ar[r] &0 \\
0 \ar[r] &\Oo_S c_0 \ar[u] \ar[r]  & \widehat{\theta}_{\Ll/S} \ar[u]^{{T}_\h} \ar[r] & \theta_{\Ll/S} \ar[u] \ar[r] &0 
}
\]
where the first vertical arrow maps $c_0$ to $c \dim(\g) (\check{h}+c)^{-1} \cdot \text{Id}$ and by abuse of notation we wrote $T_\h$ instead of $T_\h(-) \, \circ $.

\begin{rmk} \label{rmk-fockgexplicit} We give a local description of how the action looks like. Choose for this purpose a local parameter $t$ of $\mathcal{I}_\sigma$ so that we can associate to $D \in \theta_{\Ll/S}$ the element $\widehat{D} \in \widehat{\theta}_{\Ll/S}$. Let $X_r \circ \cdots \circ X_1$ be representatives of an element of $\Ff^+(\hL)$ whith $X_i \in \hL$. Then the action of $\theta_{\Ll/S}$ is described as follows
\[T_\h(\widehat{D}) \circ X_r \circ \cdots \circ X_1 = \sum_{i=1}^r X_r \circ \cdots \circ D(X_i) \circ \cdots \circ X_1 + X_r \circ \cdots \circ X_1 \circ T_\h(\widehat{D})
\] where $D(X_i)$ denotes the image of $X_i$ under coefficientwise derivation by $D$ (Remark \ref{rmk-coeffder}). \end{rmk}

\subsubsection{Projective representation of $\theta_{\Ll,S}$} 	\label{subsubsec-extendcentralext}

We define in this section a map of Lie algebras $\PP T_{\h,S} \colon \theta_{L,S} \to \PP\End(\Hhn{\V})$ which is induced by $T_\h$ and which will lead, in a second time, to a projective connection on the sheaves of covacua. By duality, this will induce a projective connection on conformal blocks. The construction of $\PP T_{\h,S}$ in the classical case is the content of \cite[Corollary 3.3]{looijenga2005conformal}.

Let $F^0\theta_{\Ll/S}$ be the subsheaf of $\theta_{\Ll/S}$ given by those derivations $D$ such that $D(F^1\!\Ll) \in F^1\!\Ll$, and similarly we set $F^0\theta_{\Ll,S}$ to be the subsheaf of $\theta_{\Ll,S}$ whose elements $D$ satisfy $D(F^1\!\Ll) \in F^1\!\Ll$. 

\begin{rmk} \label{rmk-coeffcasimiroobstr} Assume that $\Ll=\Oo_S(\!(t)\!)$ so that every element of $F^0\theta_{\Ll/S}$ is written as $D=\sum_{i \geq 0} a_i t^i d/dt$. The element $T_\h(\widehat{D})$ acts on $\V$ as the operator $\dfrac{a_0}{-2(\ell+\check{h})}\mathfrak{c}$ where $\mathfrak{c}$ is the Casimir element of $\hL$, hence the action is by scalar multiplication. Combining this with Remark \ref{rmk-fockgexplicit}, we obtain that $F^0(\theta_{\Ll/S})$ acts on $\Hhn{\V}$ by coefficientwise derivation up to scalars.
\end{rmk}

As in the classical case, also in our context this observation is the key input to define $\PP T_{\h,S}$. In fact we let $F^0\theta_{\Ll,S}$ act on $\Hhn{\V}$ by coefficientwise derivation so that we obtain a map
\[F^0\theta_{\Ll,S} \times \theta_{L/S} \to \PP\End(\Hhn{\V})\]
which uniquely defines the Lie algebra homomorphism
\[\PP T_{\h,S} \colon \theta_{\Ll,S} \to \PP\End(\Hhn{\V})\]
and hence the central extension $\widehat{\theta}_{\Ll,S} \to \theta_{\Ll,S}$ and the map $T_{\h,S} \colon \widehat{\theta}_{\Ll,S} \to \End(\Hhn{\V})$.

    \proof We only have to prove that the Lie algebra generated by $F^0\theta_{\Ll,S}$ and $\theta_{L/S}$ is $\theta_{\Ll,S}$. This can be checked locally, where the choice of a local parameter $t$ allows us to split the exact sequence
    \[ 0 \to \theta_{\Ll/S} \to \theta_{\Ll,S} \to \T_{S/k} \to 0, \] hence to write $\theta_{\Ll,S}$ as $\theta_{\Ll|S} \oplus \T_{S/k}$. We can in fact decompose every element $D \in \theta_{\Ll,S}$ as $D_{ver} \oplus D_{hor}$ which are uniquely determined by the conditions
    \[(\star) \quad D_{ver} \in \theta_{\Ll|S}, D_{hor}(t) = 0 \AND D_{hor}(s)=D(s) \text{ for all }s \in \Oo_S.\quad \]
    This implies that $F^0\theta_{\Ll,S}=F^0\theta_{\Ll/S} \oplus \T_{S/k}$, concluding the argument. \endproof

\begin{rmk} \label{rmk-bracketDerive} Assume that $\Ll=R(\!(t)\!)$ so that we can write every element $D \in \theta_{\Ll,S}$ as $D=D_{ver}+D_{hor}$ satisfying $(\star)$. Then Remark \ref{rmk-fockgexplicit} tells us that the action of $D$ on $\Hhn{\V}$ is given by componentwise derivation by $D$ plus right multiplication by $T(\widehat{D}_{ver})$. \end{rmk}

\begin{rmk} We want to remark that in the case in which $\V$ is the trivial representation, then the central extension $\widehat{\theta}_{\Ll,S}$ is isomorphic to $\widehat{\theta}_{\Ll/S} \oplus \T_{S/k}$, viewed as a Lie subalgebra of $\gl(\Hhne)$, where the action of $\T_{S/k}$ is by coefficientwise derivation. In fact in the previous proof we saw that locally on $S$, and up to the choice of a local parameter this is the case. By looking at Remark \ref{rmk-coeffcasimiroobstr} and the previous proof, we note that the obstruction to deduce this statement globally lies in the action of the Casimir element $c_\h$ on $\V$. When $\V$ is the trivial representation $c_\h$ acts as multiplication by zero, hence there is no obstruction. In particular, the central charge $c_0=(0, -\hslash) \in \widehat{\theta}_{\Ll,S}$ acts by multiplication by $\dim(\g) \ell/ (\ell+\check{h})$.
\end{rmk}

\subsection{The projective connection on $\VV_\ell(\V)$}    \label{subse-projconn}

The aim of this section is to induce, from $\PP T_{\h,S}$, the projectively flat connection $\nabla \colon \T_{S/k}(-\log(\Delta)) \to \PP \End(\VV_\ell(\V)_X)$. In Proposition \ref{prop-exseqALog} we realised $\T_{S/k}(-\log(\Delta))$ as the quotient $\theta_{\Aa,S}(-\Rr) / \theta_{\Aa/S}(-\Rr)$, so that the content of this section can be summarized in the following statement.

\begin{lem} \label{thm-connect} The actions of $\theta_{\Aa,S}$ and of $\theta_{\Aa/S}(-\Rr)$ on $\Hhn{\V}$ induce a projective action of $\T_S(-\log(\Delta))$ on $\VV_\ell(0)_X$. In particular $\VV_\ell(\V)_X$ is locally free if restricted to $S \setminus \Delta$.
\end{lem} 

As a consequence of it, we obtain that $\VV_\ell(\V)$ is locally free on $\Hur$. This is the first step to the proof of Theorem \ref{thm-LocFree}.

\begin{cor} \label{cor-ConnectionLocFree} The sheaf $\VV_\ell(\V)$ on $\bHur$ is equipped with a projectively flat connection with logarithmic singularities along the boundary $\Delta_{\\univ}$. In particular $\VV_\ell(\V)$ is locally free on $\Hur$.
\end{cor}

\proof As pointed out in Subsection \ref{subsec-tangHur} the tangent space of $\bHur$ at a versal covering $(\Xt \overset{q}\to X \overset{\pi}\to S, \sigma)$ is identified via the Kodaira-Spencer map with the tangent bundle $\T_{S/k}(-\log(\Delta))$. The previous theorem gives the projective action of the latter on $\VV_\ell(\V)_X$, concluding in this way the argument.
\endproof

\subsubsection*{Proof of Lemma \ref{thm-connect}} We first of all prove that the action of $\theta_{\Aa,S}(-\Rr)$ on $\Hhn{\V}$ descends to $\VV_\ell(\V)_X$. 

\begin{prop} The projective action of $\theta_{\Aa,S}(-\Rr)$ on $\Hhn{\V}$ preserves the space $\hA \circ \Hhn{\V}$, hence induces a projective action on $\VV_\ell(\V)_X$.
\end{prop}

\proof By the local description of the action of $\theta_{\Aa,S}(-\Rr)$ explained in Remark \ref{rmk-bracketDerive}, it suffices to show that the action of $\theta_{\Aa,S}(-\Rr)$  on $\hA$ by coefficientwise derivation is well defined. This follows from Remark \ref{rmk-coeffder}.
\endproof

We denote by $\PP T_{\hA,S}$ the morphism $\theta_{\Aa,S}(-\Rr) \to  \VV_\ell(\V)_X$ induced by $\PP T_{\h,S}$. To conclude the proof of lemma \ref{thm-connect} we are left to show the following proposition. 

\begin{prop} \label{prop-aquot} The morphism $\PP T_{\hA,S} \colon \theta_{\Aa,S}(-\Rr) \to \PP \End(\VV_\ell(\V)_X)$ factorizes through \[\PP T_{\hA,S} \colon \T_S(-\log(\Delta)) = \theta_{\Aa,S}(-\Rr)/\theta_{\Aa/S}(-\Rr) \to \PP \End(\VV_\ell(\V)_X).\]
\end{prop}

\proof We need to prove that $\theta_{\Aa/S}(-\Rr)$ acts on $\VV_\ell(\V)_X$ by scalar multiplication. As this can be checked locally, we can assume to have a local parameter, so that we can associate to $D \in \theta_{\Aa/S}$ the element $\widehat{D} \in \widehat{\theta}_{\Ll/S}$. We need to prove that, up to scalars, $T_{\h}(\widehat{D})$ lies in the closure of $\hA \circ \hL$ in $\overline{U}\HhL[(c+\check{h})^{-1}]$. 

For this purpose we use the description of $\widehat{C}_{\h}(\widehat{D})$ provided in Remark \ref{rmk-explLoc}. Let consider the orthonormal bases with respect to $\Res_{(\,|\,)}$ given by elements $\{\alpha_i, \beta_j\}$ and $\{a_i, b_j\}$ of $\omega_{\Ll/S} \otimes \hL$ and $\hL$ and with $a_i \in \hA$. From Remark \ref{rmk-explLoc} we can write 
\[T_\h(\widehat{D}) = \sum \circdot \, D(\alpha_i) \circ a_i \, \circdot + \sum \circdot \, D(\beta_j) \circ b_j \,\circdot.
\]
Observe that up to an element in $c\Oo_S$ we have the equality $\sum \circdot \, D(\alpha_i) \circ a_i \, \circdot = \sum a_i \circ D(\alpha_i)$, so that to conclude it is enough to show that $D(\beta_j) \in \hA$. To do this, we first need to identify where $\beta_j$'s live. Since the basis is $\Res_{(\,|\,)}$-orthonormal we know from Lemma \ref{lem-annullatore} that $(1 \otimes(\,|\,)_\hL)(\beta_j) \in \omega_\Aa \otimes \hA^\vee$. 

Recall that in Lemma \ref{lem-decomp}, we decomposed $\h$ as $\bigoplus_{i=0}^{p-1} g^{(-i)}\otimes \E_{i}$. Using this decomposition, and the fact that $(\,|\,)_\h$ provides an isomorphism between $\g^{(-i)}\otimes_k \E_i$ and $(\g^{(i)} \otimes_k \E_{p-i}(\Rr))^\vee$ for $i \neq 0$, we deduce that 
\[\beta_j \in \left(\g^\Gamma \pi_*\Oo_{X^*} \oplus \bigoplus_{i=1}^{p-1} \left(\g^{(i)} \otimes_k \pi_* \E_{p-i}(\Rr)|_{X^*}\right) \right)\otimes \omega_\Aa\]
It follows that 
\[D(\beta_j) \in \g^\Gamma \pi_*\Oo_{X^*}(-\Rr) \oplus \bigoplus_{i=1}^{p-1} (\g^{(i)} \otimes_k \pi_* \E_{p-i}|_{X^*})\]
and hence is contained in $\hA= \underset{i=0}{\overset{p-1}{\bigoplus}} \g^{(-i)} \otimes_k \pi_*\E_{i}|_{X^*}$.\endproof

\subsubsection{Projective connection on $\VV_\ell(\V_1, \dots, \V_n)$} \label{subsec-SemiLocalCase2}

As we did in Section \ref{subsec-SemiLocalCase}, whose notation we are going to use here, we extend also the projective connection to the case in which more points on $X$ are fixed. Observe first of all that the identification of the tangent space of $\bHur[n]$ at a versal covering $(\Xt \overset{q}\to X \overset{\pi}\to S, \{\sigma_i\})$ with $\T_{S/k}(-\log(\Delta))$ still holds. Since Proposition \ref{prop-exseqALog} still holds, this implies that we are allowed to provide the projective connection on $\VV_\ell(\V_1, \dots, \V_n)_X$ in terms of a projective action of $\theta_{\Aa,S}(-\Rr)$ on $\VV_\ell(\V_1, \dots, \V_n)_X$.

We denote by $\theta_{\Ll/S}$ the direct sum of $\theta_{\Ll_i/S}$, and we obtain a central extension $\widehat{\theta}_{\Ll/S}$ thereof as the quotient of the direct sum of $\widehat{\theta}_{\Ll_i/S}$ which identifies $(0,\hslash_i) \in \widehat{\theta}_{\Ll_i/S}$ with $(0,\hslash_j) \in \widehat{\theta}_{\Ll_j/S}$.
The Sugawara representation $T_\h \colon \widehat{\theta}_{\Ll/S} \to \overline{U}\HhL[(c+\check{h})^{-1}]$ is induced from the Sugawara representations of $\widehat{\theta}_{\Ll_i/S}$ and gives the projective action of $\theta_{\Ll,S}$ on $\Hhn{\V_1, \dots, \V_n}$.

Combining all these elements with the case $n=1$ we obtain the following generalization of Corollary \ref{cor-ConnectionLocFree}.

\begin{cor} \label{cor-HhnLog} For every $i \in \{1, \dots, n\}$ let  $\V_i \in \Repl({\sigma_{i,un}}^*\h)$. The coherent module $\VV_\ell(\V_1, \dots, \V_n)$ is equipped with a projective action of $\T_{\bHur[n]}(-\log(\Delta))$. In particular it is locally free over $\Hur[n]$. \end{cor}

    \section{Factorization rules and propagation of vacua} \label{sec-factorizationRulesS}              %

We show in this section that the sheaf $\VV_\ell(0)$ descends to $\Hurz$ by means of the propagation of vacua, and we provide the factorization rules which compare the fibre of $\VV_\ell(0)$ over a nodal curve $X$ with the fibres of the sheaves $\VV_\ell(\V)$ on its normalization $X_N$. We will proceed following the approach of \cite[Section 4]{looijenga2005conformal}. 

\subsection{Independence of number of sections}

In this section we want to show that the sheaf $\VV_\ell(0)$ actually descends to a vector bundle on $\Hurz$ as a consequence of Proposition \ref{prop-PropVacua}. Following \cite[Proposition 2.3]{beauville1994conformal} we state and prove the aforementioned proposition, called also \emph{propagation of vacua}, because it shows that we can modify the sheaf of covacua by adding as many sections as we want to which we attach the trivial representation to obtain a sheaf isomorphic to the one we started with. 

\begin{setting} In this section we fix the following objects.
\begin{itemize}
\item Let $(\Xt \overset{q}{\to} X \overset{\pi}{\to} S=\Spec(R),\sigma_1, \dots, \sigma_n, \sigma_{n+1}, \dots, \sigma_{n+m})$ be an element of $\Hur[n+m](S)$. 
\item Denote by $\B:=\Oo_{X\setminus \{S_1, \dots, S_n\}}$ and by $\A:=\Oo_{X\setminus \{S_1, \dots, S_n, S_{n+1}, \dots, S_{n+m}\}}$ and set $\hB:=\pi_*(\h \otimes \B)$ which is contained in $\hA:=\pi_*(\h\otimes \A)$.
\item For every $i \in \{1, \dots, n\}$ fix $\V_i \in \Repl(i):=\Repl(\sigma_i^*\h)$ and for every $j \in \{1, \dots, m\}$ we fix $\W_j \in \Repl(n+j):=\Repl(\sigma_{n+j}^*\h)$. 
\end{itemize}\end{setting}

Under these conditions we notice that $\hB$ acts on each $\W_j$ since ${\sigma_{n+j}}^*\pi^*\hB$ maps naturally to ${\sigma_{{n+j}}}^*\h$ and the latter acts on $\W_j$ by definition. 

\begin{prop}\label{prop-PropVacua} The inclusions $\W_j \to \Hhn{\W_j}$ induce an isomorphism
\[\hB \Qsx \left(\Hhn{\V_1, \dots, \V_n} \otimes \bigotimes_{j=1}^m \W_j\right) \cong \hA \Qsx \Hhn{\V_1, \dots, \V_n, \W_1, \dots, \W_m}
\]of $\Oo_S$-modules.
\end{prop}

\proof We sketch here the main ideas of the proof, using the same techniques used by Beauville in \cite[Proof of Proposition 2.3]{beauville1994conformal}. By induction it is enough to prove the assertion for $m=1$, so that we need to prove that the inclusion $\W \to \Hhn{\W}$ induces an isomorphism 
\[\phi \colon \hB \Qsx \left(\Hhn{\V_1, \dots, \V_n} \otimes \W \right) \overset{\cong}\longrightarrow \hA \Qsx \Hhn{\V_1, \dots, \V_n, \W}.\]
The morphism is well defined on the quotients as the inclusion of $\hB$ in $\hL_{n+1}$ factors through $\hB \to \hA$. Since the inclusion $\W \to \Hhn{\W}$ factors through $\tHh{\W}$, we prove the proposition in two steps.

\textit{\textbf{Claim 1.} The inclusion $\W \to \tHh{\W}$ induces an isomorphism}
\[\widetilde{\phi} \colon \hB \Qsx \left(\Hhn{\V_1, \dots, \V_n} \otimes \W \right) \longrightarrow \hA \Qsx \left( \Hhn{\V_1, \dots, \V_n} \otimes \tHh{\W}\right).\]

\textit{\textbf{Claim 2.} The projection map 
\[\hA \Qsx \left( \Hhn{\V_1, \dots, \V_n} \otimes \tHh{\W}\right) \longrightarrow
\hA \Qsx  \Hhn{\V_1, \dots, \V_n,\W}\]
is an isomorphism.}

We give the proof of \textbf{Claim 1}, as for \textbf{Claim 2} one can refer to \cite[(3.4)]{beauville1994conformal}. We just remark
that in the proof of Claim 2 it is used that the level of $\W$ is bounded by $\ell$ and the local description of the maximal submodule $\Z$ of $\tHh{\W}$.

Since checking that $\widetilde{\phi}$ is an isomorphism can be done locally on $S$, there is no loss in generality in assuming that the $\mathcal{I}_{i}$'s are principal so that $\Op \cong \bigoplus R[\![t_i]\!]$ and that there are isomorphisms $\hLi \cong \gL_i$. Observe that the exact sequence of $R$-modules
\[ 0 \to \hB \to \hA \to \hA/ \hB \to 0
\]
splits because the quotient $\h^-:=\hA/ \hB$ is isomorphic to $\hLi[n+1] / \h_{\Op_{n+1}}$ which can be identified with $\g \otimes_k R[t_{n+1}^{-1}]t^{-1}_{n+1}$. Moreover we observe that since the residue pairing is trivial on $\h^-$, this is a Lie sub algebra of $\hA$, hence we are left to prove that \[\Hhn{\V_1, \dots, \V_n} \otimes \W \longrightarrow \h^- \Qsx \left( \Hhn{\V_1, \dots, \V_n} \otimes \tHh{\W}\right)\]
is an isomorphism. Observe that this statement no longer depends on the covering $\Xt \to X$, so once we choose isomorphisms $\hLi \cong \gL_i$, this follows from the classical case.\endproof 

As previously announced, this proposition has important consequences.

\begin{cor} \label{cor-PropVacua} For all $n$ and $m \in \ZZ_{> 0}$ there is a natural isomorphism
\[ \left(\textbf{Forg}_{n+m,n}\right)^*\left(\VV_\ell(\V_1, \dots, \V_n)\right) \cong  \VV_\ell(\V_1, \dots, \V_n, 0, \dots, 0)
\] of vector bundles on $\Hur[n+m]$.\qed
\end{cor}

Which leads to the following result.

\begin{cor} \label{cor-descent} The vector bundle $\VV_\ell(0)$ defined on $\Hur$ descends to a vector bundle on $\Hurz$.\end{cor}

\proof We can construct the sequence
\[\xymatrix{ \Hur[3] \ar@<-.6ex>[r] \ar[r] \ar@<+.6ex>[r] & \Hur[2] \ar@<-.3ex>[r]_{f_1} \ar@<+.3ex>[r]^{f_2} & \Hur \ar[r] &\Hurz }\]
where the horizontal morphisms, which are faithfully flat, are given by forgetting one of the sections. The vector bundle $\VV_\ell(0)$ on $\Hur$ then descends from $\Hur$ to $\Hurz$ because Corollary \ref{cor-PropVacua} provides a canonical isomorphism $\phi_{12}$ between $f_1^*\VV_\ell(0)$ and $f_2^*\VV_\ell(0)$. The compatibility of the isomorphisms $\phi_{ij}$ on $\Hur[3]$ holds by construction.
\endproof

\begin{rmk} \label{rmk-dropaffine} When we defined $\VV_\ell(\V)$ on $\bHur$, we assumed that the complement of the section $\sigma$ was affine. Corollary \ref{cor-PropVacua} allows us to remove this assumption: in fact if this is not the case, we can add finitely many sections to which we attach the trivial representation and set $\VV_\ell(\V)$ to be $\VV_\ell(\V, 0, \dots,0)$. The same holds for $\VV_\ell(\V_1, \dots, \V_n)$ on $\bHur[n]$.
\end{rmk}

\subsection{Nodal degeneration and factorization rules} \label{subsec-NodDeg}

In this section we want to compare the sheaf of covacua $\VV_\ell(0)_X$ attached to a covering of nodal curves $\Xt \to X$, to the sheaves of the form $\VV_\ell(\V)_{X_N}$, attached to the normalization $\Xt_N \to X_N$ of the covering we started with.

\begin{setting} We will consider the following objects.
\begin{itemize}
\item Let $(\Xt \overset{q}\to X \overset{\pi}\to \Spec(k), \mathfrak{p}_1, \dots, \mathfrak{p}_n) \in \bHur[n](\Spec(k))$ and assume that $X$ has one double point $\p \in X(k)$. 
\item Let $X_N$ be the partial normalization of $X$ separating the branches at $x$ and set $q_N \colon \Xt_N:=\Xt \times_X X_N \to X_N$. The points of $X_N$ mapping to $\p$ are denoted $\p_+$ and $\p_-$.

\end{itemize} \end{setting}

\begin{rmk} Observe that $q_N$ is a $\Gamma$-covering with the action of $\Gamma$ induced by the one on $\Xt$, so that it is ramified only over $\Rr$. The Lie algebra $\h_N := \h \times_X X_N$ is then isomorphic to the Lie algebra of $\Gamma$-invariants of ${q_N}_*(\g \otimes_k \Oo_{\Xt_N})$. Furthermore, the normalization provides an isomorphism between the $k$-Lie algebra $\h|_{\p}$ and $\h_N|_{\p_\pm}$.
\end{rmk}

Let $\pi_N \colon X_N \to \Spec(k)$ be the structural morphism and consider the Lie algebras \[\h_{\Aa_N}:={\pi_N}_*\left(\h_N \otimes \Oo_{X_N \setminus \{\p_1, \dots, \p_n\}}\right)\] and \[\h_{\Ll_N} := \bigoplus_{i=1}^n \varinjlim_{M \in \ZZ_{\geq 0}} {\pi_N}_*\varprojlim_{m \in \ZZ_{> 0}} \h_N \otimes \mathcal{I}_i^{-M}/ \mathcal{I}_i^m\] which are the analogues of $\hA$ and $\hL$ for the marked covering $(\Xt_N \to X_N, \{\p_i\})$. Observe that since $X$ and $X_N$ are isomorphic outside of $\p$, the Lie algebras $\h_{\Ll_N}$ and $\hL$ are isomorphic.

As observed in the previous remark, since $\h|_{\p}$ and $\h_N|_{\p_\pm}$ are isomorphic, every representation $\W$ of $\h|_\p$, is also a representation of $\h_N|_{\p_\pm}$. Let denote by $\W^*$ the dual of $\W$ and view $\W \otimes_k \W^*$ as a representation of $\h_N|_{\p_+} \oplus \h_N|_{\p_-}$, with $\h_N|_{\p_+}$ acting on $\W$ and $\h_N|_{\p_-}$ on $\W^*$. This induces an action of $\h_{\Aa_N}$ on $\W \otimes_k \W^*$ as
\[\alpha *(w\otimes \phi)=[X]_{\p_+}w \otimes \phi + w \otimes [X]_{\p_-}\phi\] where $[\star]_{\p_\pm}$ denotes the reduction modulo the ideal defining $\p_\pm$. Let $b_\W$ denote the trace morphism $\End(\W)=\W \otimes_{\Oo_S} \W^* \to \Oo_S$ which is compatible with the action of $\h|_\p$. We can formulate the \emph{factorization rules} controlling the nodal degeneration as follows.

\begin{prop} \label{prop-NodalDegenerationH}  The morphisms $\{b_\W\}$ induce an isomorphism \[ \bigoplus_{\W \in \Repl(\h|_\p)} \h_{\Aa_N} \Qsx \left(\Hhn{\V_1, \dots, \V_n} \otimes (\W \otimes \W^*) \right) \to \hA \Qsx \Hhn{\V_1, \dots, \V_n}\]\end{prop}

The proof of this result is a mild generalization of the proof of \cite[Proposition 6.1]{looijenga2005conformal}, which in turn is a consequence of Schur's Lemma. We give an overview of it.

\proof Fix an isomorphism between $\h|_\p$ and $\g$ so that $\Repl(\h|_{\p})$ is identified with $P_\ell$. Denote by $\Spec(\Aa)=X\setminus\{\p_1, \dots, \p_n \}$ and $\Spec(\Aa_N)=X_N \setminus \{\p_1, \dots, \p_n \}$ and let $I_\p \subset \Aa$ be the ideal defining $\p$, so that the normalization gives the diagram
\[\xymatrix{
0 \ar[r]& I_\p \ar[r] \ar@{=}[d]& \Aa \ar[r]\ar[d]& k \ar[r] \ar[d]^\Delta& 0\\
0 \ar[r]& I_\p \ar[r] & \Aa_N \ar[r]& k\oplus k \ar[r] & 0
}\]
whose rows are exact. As in the classical case, we consider a similar diagram of Lie algebras, which, in contrast to the classical case, is not obtained by tensoring by $\g$. Define $\h_I$ as the tensor product $I_\p \otimes_{\Aa} \hA$, and observe that the quotient $\hA/\h_{I_\p}$ is $\h|_{\p}$, which is then isomorphic to $\g$. Repeating the construction on $X_N$, we obtain the commutative diagram of Lie algebras
\[\xymatrix{
0 \ar[r]& \h_{I_\p} \ar[r] \ar@{=}[d]& \hA \ar[r]\ar[d]& \g \ar[r] \ar[d]^\Delta& 0\\
0 \ar[r]& \h_{I_\p} \ar[r] & \h_{\Aa_N} \ar[r]& \g\oplus \g \ar[r] & 0.
}\]

Consider the $\h_{\Aa_N}$-module $\Hhn{\V_1, \dots, \V_n}$ and observe that its quotient $M:= \h_{I_\p} \Qsx \Hhn{\V_1, \dots, \V_n}$ is a finite dimensional representation of $\g \oplus \g$, because the quotient $\hA \Qsx \Hhn{\V_1, \dots, \V_n}$ is finite dimensional and the quotient $\hA / \h_{I_\p}$ is one dimensional. It is moreover a representation of $\g \oplus \g$ of level less or equal to $\ell$ relative to each factors. Since $\h_{I_\p}$ acts trivially on $\W \otimes \W^*$, the maps $\{b_\W\}$ induce the morphism 
\[ \bigoplus_{W \in P_\ell} \h_I \Qsx \Hhn{\V_1, \dots, \V_n} \otimes (\W \otimes \W^*) \to \h_I \Qsx \Hhn{\V_1, \dots, \V_n}\]
Observe, using the above diagram, that if we consider $M$ as a $\g$-module via the diagonal action, and we denote this $\g$-representation by $M^\Delta$, then $\g \Qsx M^\Delta$ is exactly $\hA \Qsx \Hhn{\V_1, \dots, \V_n}$. After these considerations, the proof of the proposition boils down to showing that if $M$ is a finite dimensional representation of $\g \oplus \g$ of level at most $\ell$, then the morphisms $\{b_\W\}$ induce the isomorphism 
\[\bigoplus_{W \in P_\ell} \g \oplus \g \Qsx (M \otimes (W \otimes W^*)) \to \g \Qsx M^\Delta.\]
Schur's Lemma ensures that the set of morphisms between irreducible Lie algebra representations is a skew field, and since without loss of generality we might assume $M$ to be an irreducible $\g \oplus \g$ representation of the form $V_1 \otimes V_2$ for $V_i \in P_\ell$, we conclude. \endproof

Denote by $\h_{{\Aa_N}^*}$ the Lie algebra ${\pi_N}_*(\h_N \otimes \Oo_{X_N \setminus \{\p_1, \dots, \p_n, \p_+, \p_-\}})$. Then Proposition \ref{prop-PropVacua} allows us to rewrite the previous proposition as an isomorphism
\[ \bigoplus_{\W \in \Repl(\h|_\p)} \h_{{\Aa_N}^*} \Qsx \Hhn{\V_1, \dots, \V_n, \W, \W^*} \to \hA \Qsx \Hhn{\V_1, \dots, \V_n}\]
and in particular implies the isomorphism \[ \bigoplus_{\W \in \Repl(\h|_\p)} \VV_\ell(\W,\W^*)_{X_N} \to \VV_\ell(0)_X\]
where we see $X_N$ naturally marked by $\p_+$ and $\p_-$.

\begin{rmk} Observe that in the case in which $x$ is a non separating node, i.e. when $X_N$ has the same number of irreducible components of $X$, then the genus of $X_N$ is one less the genus of $X$. By repeatedly applying the factorization rules, it is possible then to reduce the computation of the rank of sheaves of covacua on $\bHur[n]$ to the case of covacua over $\overline{\mathcal{H}\textrm{ur}}(\Gamma, \xi)_{0,N}$.
\end{rmk}

\begin{rmk} Let $(\Xt \to X \to S, \sigma) \in \bHur(S)$ and assume that it is possible to normalize the family (for example assuming that the nodes of $X$ are given by a section $\varsigma \colon S \to X$). Then Proposition \ref{prop-NodalDegenerationH} still holds by replacing the index set $\Repl(\h|_\p)$ with $\Repl(\varsigma^*\h)$.
\end{rmk}

\section{Locally freeness of the sheaf of conformal blocks} \label{sec-locfree}

In this section we show how to use a refined version of the factorization rules to prove that the sheaves of covacua $\VV_\ell(\V_1, \dots, \V_n)$, and hence the conformal blocks $\VV_\ell(\V_1, \dots, \V_n)^\dagger$, are locally free also on the boundary of $\bHur[n]$, concluding in this way the proof of Theorem \ref{thm-LocFree}. For simplicity only we will assume $n=1$.

\subsection{Canonical smoothing} \label{subsec-Canonicalsmoothing}

As previously stated, we want to prove that $\VV_\ell(\V)$ is a locally free sheaf on $\bHur$. For this purpose we describe here a procedure to realise a covering of nodal curves as the special fibre of a family of coverings which is generically smooth. The idea is to induce a deformation of the covering $\Xt \to X$ from the canonical smoothing of the base curve $X$ provided in \cite{looijenga2005conformal}. As already noted in Remark \ref{rmk-NCD}, it is essential that the branch locus $\Rr$ of the covering $q \colon \Xt \to X$ is contained in the smooth locus of $X$.

Let $(\Xt \overset{q_0}\to X \overset{\pi_0}\to \Spec(k),\sigma_0) \in \bHur(\Spec(k))$ with $\p\in X(k)$ the unique nodal point of $X$. The goal of this section is to construct a family $\Xxt \to \Xx$ belonging to $\bHur(\Spec(k[\![\tau]\!]))$ which deforms $(\Xt \to X)$ and whose generic fibre is smooth, i.e. it lies in $\Hur(\Spec(k(\!(\tau)\!)))$. 

\subsubsection{The intuitive idea} \label{subsubsec-intuito} The idea which is explained in \cite{looijenga2005conformal}, is to find a deformation $\Xx$ of $X$ which replaces the formal neighbourhood $k[\![\tp,\tm]\!]/\tp\tm$ of the nodal point $\p$ with the $k[\![\tau]\!]$-algebra $k[\![\tp,\tm, \tau]\!] /\tp\tm=\tau$. This can be achieved with the following geometric construction. We first normalize the curve $X$ obtaining the curve $X_N$ with two points $\p_+$ and $\p_-$ above $\p$. We blow up the trivial deformation $X_N[\![\tau]\!]$ of $X_N$ at the points $\p_+$ and $\p_-$ and note that the formal coordinate rings at $\p_\pm$ in the strict transform are of the form $k[\![\tpm, \tau/\tpm]\!]$. We then obtain the neighbourhood $k[\![\tp,\tm, \tau]\!] /\tp\tm=\tau$ by identifying $\tp$ with $\tau/\tm$. The deformation $\Xxt$ of $\Xt$ is induced from the one of $X$ because the singular point $\p$ does not lie in $\Rr$.

\subsubsection{Construction of $\Xxt \to \Xx$} We will realise the canonical smoothing of $\Xt \to X$ it by constructing compatible families \[(\Xxt^n \overset{q_n}\to \Xx^n \overset{\pi_n}\to \Spec(k[\tau]_n) ,\sigma_n) \in \bHur(\Spec(k[\tau]_n))\] where $k[\tau]_{n} := k[\tau]/(\tau^{n+1})$ for $n \in \ZZ_{\geq 0}$. As these are infinitesimal deformations, we only need to change the structure sheaf. As we have previously done, we normalize both $\Xt$ and $X$ obtaining $(\Xt_N \overset{q} \to  X_N \overset{\pi} \to \Spec(k),\sigma_0,\p_+,\p_-) \in \Hur[3](\Spec(k))$. We fix furthermore local coordinates $\tp$ and $\tm$ at the points $\p_+$ and $\p_-$. 

Let $U$ be an open subset of $X$ and $n \in \ZZ_{\geq 0}$. If $U$ does not contain $\p$ we set $\Oo_{\Xx^n}(U):=\Oo_X(U)[\tau]/\tau^{n+1}$. Otherwise, if $\p \in U$, we set
\[\xymatrix{\Oo_{\Xx^n}(U):=\ker \left( \dfrac{k[\![\tp,\tm]\!][\tau]}{\tp\tm=\tau, \tau^{n+1}} \oplus \Oo_{\Xx^n}(U\setminus \{ \p \}) \ar[rr]^-{\alpha_n-\beta_n}\right. && \left.\dfrac{k(\!(\tp))\![\tau]}{\tau^{n+1}} \oplus \dfrac{k(\!(\tm)\!)[\tau]}{\tau^{n+1}}\right)
}\]
where 
\[ \alpha_n \colon \dfrac{k[\![\tp,\tm]\!][\tau]}{\tp\tm=\tau, \tau^{n+1}} \longrightarrow \dfrac{k(\!(\tp)\!)[\tau]}{\tau^{n+1}} \oplus \dfrac{k(\!(\tm)\!)[\tau]}{\tau^{n+1}}\] is the $k[\tau]_n$-linear morphism given by $\tp \mapsto (\tp, (\tm)^{-1}\tau)$ and $\tm \mapsto ((\tp)^{-1}\tau, \tm )$, and 
\[\beta_n \colon \Oo_{\Xx^n}(U\setminus \{ \p \}) \longrightarrow \dfrac{k(\!(\tp)\!)[\tau]}{\tau^{n+1}} \oplus \dfrac{k(\!(\tm)\!)[\tau]}{\tau^{n+1}}
\] 
sends $\psi \in \Oo_{\Xx^n}(U\setminus \{ \p \})$ to $(\psi_{+}, \psi_{-})$ where $\psi_{\pm}$ is the expansion of $\psi$ at the point $\p_\pm$ using the identifications $\Oo_{\Xx^n}(U\setminus \{ \p \})=\Oo_{X}(U\setminus \{ \p \})[\tau]/\tau^{n+1}=\Oo_{X_N}(U_N\setminus \{ \p_+, \p_-\})[\tau]/\tau^{n+1}$.

\begin{rmk} \label{rmk-neightau} Observe that the completion of $\Oo_{\Xx^n}$ at the point $\p$ is isomorphic to $k[\![\tp,\tm]\!][\tau]/(\tp\tm=\tau, \tau^{n+1}) = k[\![\tp,\tm]\!]/(\tp\tm)^{n+1}$. In fact note that once we take the completion of $\Oo_{\Xx^n}(U\setminus \p)$ at the point $\p$ we obtain exactly $\dfrac{k(\!(\tp)\!)[\tau]}{\tau^{n+1}} \oplus \dfrac{k(\!(\tm)\!)[\tau]}{\tau^{n+1}}$, the map $\beta_n$ becoming the identity. The kernel of $\alpha_n-\beta_n$ is then identified with $k[\![\tp,\tm]\!][\tau]/(\tp\tm=\tau, \tau^{n+1})$ as claimed.

Observe furthermore that once we take the limit for $n \to \infty$, then the formal neighbourhood of $\p$ will be $k[\![\tp,\tm,\tau]\!]/\tp\tm=\tau$ as asserted in the subsection \ref{subsubsec-intuito}. The map $\alpha_n$ describes the process of glueing the formal charts around $\p_+$ and $\p_-$.
\end{rmk}

Note moreover that for all $n \in \ZZ_{> 0}$ there are natural maps $g^n \colon \Xx^{n-1} \to \Xx^{n}$ induced by the identity on topological spaces and by the projection $k[\tau]_{n} \to k[\tau]_{n-1}$ on the structure sheaves. 

\begin{lem} For every $n \in \ZZ_{> 0}$ the family $\Xx^n$ is a curve over $\Spec(k[\tau]_n)$ deforming $X$.
\end{lem}
\proof We need to prove that $\Oo_{\Xx^n}$ is flat and proper over $\Spec(k[\tau]_n)$. Once we show that the $\Xx^n$ is of finite type, we can use the valuative criterion to deduce that $\Xx^n$ is proper over $k[\tau]_n$. Observe that the kernel of the map ${g^n}^* \colon \Oo_{\Xx^n} \to \Oo_{\Xx^{n-1}} \to 0$ is $\tau^n\Oo_{X}$. Outside $\p$, as the deformation is trivial, this is true. On an open $U$ containing $\p$, the snake lemma tells us that this is the kernel of
\[\xymatrix{\tau^n\dfrac{k[\![\tp,\tm]\!]}{\tp\tm} \oplus \tau^n\Oo_{X}(U\setminus \{ \p \}) \ar[rr]^-{\alpha_n-\beta_n} && \tau^n k(\!(\tp)\!) \oplus \tau^n k(\!(\tm)\!)}\]
where $\alpha_n$ and $\beta_n$ are the gluing functions defining $\tau^n\Oo_X$. We can conclude that $\Oo_{\Xx^n}$ is of finite type by using induction on $n$ and observing that $\Xx^0=X$. This moreover shows the flatness of the family.
\endproof

This deformation of $X$ induces a deformation of $\Xt$ which, in rough terms, is obtained as the trivial deformation outside the points $\p_1, \dots, \p_p$ of $\Xt$ lying above $\p$, and for every $j \in \{1, \dots, p\}$, the formal neigbourhood of $\Xxt_n$ around $\p_i$ will be isomorphic to the formal neighbourhood of $\Xx_n$ around $\p$.

To do this, let denote by $\p_{j,+}$ and $\p_{j,-}$ the two points of $\Xt_N$ mapping to $\p_j$. We fix local coordinates $\tpj$ at $\p_{j,\pm}$ so that \[\Xt_N \times_{X_N} \Spec(k[\![\tpm]\!]) \cong \Spec\left(\oplus_{j=1}^p k[\![\tpmj]\!]\right)\] 
and let consider $U$ an open subset of $\Xt$. If $U$ is disjoint from $q^{-1}(\p)=\{\p_1, \dots, \p_p\}$, we set $\Oo_{\Xxt^n}(U)=\Oo_{\Xt}(U)[\tau]/\tau^{n+1}$. Let $i \in \{1, \dots, p\}$ and if $U$ contains $\p_i$ but not $\p_j$ for all $j \neq i$ we set
\[\xymatrix{\Oo_{\Xxt^n}(U):=\ker \left( \dfrac{k[\![\tpj[i],\tmj[i]\!]\!][\tau]}{\tpj[i]\tmj[i]=\tau, \tau^{n+1}} \oplus \Oo_{\Xxt^n}(U\setminus \{ \p_i \}) \ar[rr]^-{\widetilde{\alpha}_n-\widetilde{\beta}_n}\right. && \left.\dfrac{k(\!(\tpj[i])\!)[\tau]}{\tau^{n+1}} \oplus \dfrac{k(\!(\tmj[i])\!)[\tau]}{\tau^{n+1}}\right)
}\]
where the maps $\widetilde{\alpha}_n$ and $\widetilde{\beta}_n$ are defined as in the case of $\Xx^n$.

\begin{rmk} As was shown for $\Xx^n$, also $\Xxt^n$ is a curve over $\Spec(k[\tau]_n)$ deforming $\Xx$. \end{rmk}

Let denote by $\Rr_n$ the trivial deformation of the branch locus $\Rr$ inside $\Xx^n$. The natural map $q_n \colon \Xxt^n \to \Xx^n$ which extends $q_0 \colon \Xt \to X$ realizes $\Xxt^n \to \Xx^n$ as a $\Gamma$-covering which is étale exactly outside $\Rr_n$ since the map $q_n$ is étale on $\p$ by construction. Furthermore, as $\sigma_0$ is disjoint from the singular locus, it follows that for every $n \in \ZZ_{\geq 0}$ we can set $\sigma_n$ to be the trivial deformation of $\sigma_0$. 

By taking the direct limit of this family of deformations we obtain the $\Gamma$-covering of formal schemes $\Xxt^\infty \to \Xx^\infty$ over $\text{Spf}(k[\![t]\!])$. To prove that $\Xxt^\infty \to \Xx^\infty$ is algebraizable, i.e. that is comes from an algebraic object $\Xxt \to \Xx \to \Spec(k[\![\tau]\!]\!])$, we can invoke Grothendieck's existence theorem (\cite[Théorème 5.4.5]{EGAIII2}) so that we are left to prove that the family $(\Xxt^n \to \Xx^n)_n$ is equipped with a compatible family of very ample line bundles. This is true because given a smooth point $P$ of $X$ which is not in $\Rr$ and $m$ sufficiently big, we know that $\Oo(mP)$ is a very ample line bundle on $X$ whose pullback to $\Oo_{\Xt}$ is also very ample. Since $P$ lies in the smooth locus of $X$ these line bundles extend naturally to very ample line bundles on $\Xx^n$ and on $\Xxt^n$, providing the wanted family of very ample line bundles.

We refer to the covering $(q \colon \Xxt \to \Xx, \sigma)$ that we have just constructed as the \emph{canonical smoothing} of $(q_0 \colon \Xt \to X, \sigma_0)$. Observe that the generic fibre of $\Xxt \to \Xx$ is a covering of smooth curves over $k(\!(\tau)\!)$ because, as one can deduce from Remark \ref{rmk-neightau}, the formal neighbourhood of $\p$ is given by $k[\![\tp,\tm]\!](\!(\tau)\!)/\tp\tm=\tau$.

\subsection{Local freeness} \label{subsec-localfreeness}

The aim of this section is to show that, in the setting of the previous section, $\VV_\ell(\V)_\Xx$ is a locally free $k[\![\tau]\!]$-module. 
We can depict the situation that we described in the previous section in the following diagram 
\[\xymatrix{
\Xt=\Xxt_0 \ar[r] \ar[d]^{q_0} & \Xxt \ar[d]^q & \ar[l] \Xxt_\eta \ar[d]^{q_\eta}\\
X=\Xx_0 \ar[r] \ar[d]^{\pi_0} & \Xx \ar[d]^\pi & \ar[l] \Xx_\eta \ar[d]^{\pi_\eta}\\
\Spec(k) \ar@/^1.2pc/[u]^{\sigma_0} \ar[r]^0 &  \Spec(k[\![\tau]\!]) \ar@/^1.2pc/[u]^\sigma	&\ar[l]_\eta  \Spec(k(\!(\tau)\!)) \ar@/^1.2pc/[u]^{\sigma_\eta}
}\]
where the covering $\Xxt_\eta \to \Xx_\eta$ is a $\Gamma$-covering of smooth curves and we denote by $\p$ the nodal point of $X$. Let $\V \in \Repl(\sigma)$ so that $\Hhn{\V}$ is a $k[\![\tau]\!]$-module. We use the subscript $0$ to denote the pullback along $0$, i.e. the restriction to the special fibre, so that $\V_0$ denotes the induced representation of ${\sigma_0}^*\h$.

\begin{rmk} \label{rmk-modulotau} Observe that there is a canonical injection of $k[\![\tau]\!]$-modules $\Hhn{\V} \to \Hhn{\V_0}[\![\tau]\!]$ which is an isomorphism modulo $\tau^n$ for every $n \in \ZZ_{> 0}$. Moreover we have by construction that $(\hA)_0$ is isomorphic to $\h_{\Aa_0}=\h_A$ and so $\left(\VV_\ell(\V)_\Xx\right)_0$ is isomorphic to $\VV_\ell(\V_0)_X$.\end{rmk}

The main result is that $\VV_\ell(\V)_\Xx$ is the trivial deformation of $\VV_\ell(\V_0)_X$ as stated in the following theorem.

\begin{thm} \label{thm-HhnTrivDef} There is an isomorphism  
\[ \VV_\ell(\V)_\Xx \cong \VV_\ell(\V_0)_X[\![\tau]\!]
\]
of $k[\![\tau]\!]$-modules. In particular $\VV_\ell(\V)_\Xx$ is a free $k[\![\tau]\!]$-module.
\end{thm}

\subsubsection{Proof of Theorem \ref{thm-HhnTrivDef}}

\begin{notation} In what follows we denote by $\Op_N:=k[\![\tp]\!] \oplus k[\![\tm]\!]$ the $k$-algebra which is the coordinate ring of the disjoint union of the formal neighbourhoods at the points $\p_\pm$ in $X_N$. Similarly $\Ll_N:=k(\!(\tp)\!) \oplus k(\!(\tm)\!)$ represents the disjoint union of the punctured formal neighbourhoods at the points $\p_\pm$ in $X_N$. Moreover we will write $k[\![\tp,\tm]\!]$ in place of $k[\![\tau, \tp, \tm]\!] \Qdx \tp\tm=\tau$. Recall that this is the completion of $\Oo_\Xx$ at the point $\p$.\end{notation}

\begin{lem} \label{Lem43} The canonical smoothing identifies $k[\![\tp,\tm]\!]$ with the subalgebra of $\Ll_N[\![\tau]\!]$ consisting of elements \[ \Op_{\p}:=\left\lbrace \sum_{i,j \geq 0} a_{ij}\left(\tp^{i-j}\tau^j\,,\,\tm^{j-i}\tau^i\right) \, | \, a_{ij} \in k \right\rbrace\]
via the map sending $\tp$ to $(\tp, \tm^{-1}\tau)$ and $\tm$ to $(\tp^{-1} \tau, \tm)$.\end{lem}
\proof Taking the limit of the definition of $\Oo_{\Xx^n}$ we identify the formal neighbourhood of $\Xx^\infty$ at $\p$ with 
\[\ker \left( k[\![\tp,\tm]\!] \oplus \Ll_N[\![\tau]\!] \overset{\alpha-Id}\longrightarrow \Ll_N[\![\tau]\!]\right)\]
where $\alpha(\tp)=(\tp, \tm^{-1}\tau)$ and $\alpha(\tm)=( \tp^{-1}\tau, \tm)$. 
\endproof

In view of Proposition \ref{prop-NodalDegenerationH} we identify $\VV_\ell(\V_0)_X[\![\tau]\!]$ with \[\bigoplus_{W \in \Repl(\h|_\p)} \h_{A_N} \Qsx \left(\Hhn{\V_0} \otimes (W \otimes W^*)\right)[\![\tau]\!]\] or equivalently with 
\[\bigoplus_{W \in \Repl(\h|_\p)} \h_{{A_N}^*} \Qsx \Hhn{\V_0, W, W^*}[\![\tau]\!].\]

Recall that $\HhL$ is a filtered Lie algebra, hence this induces a filtration on $U\HhL$ and by consequence on $\Ff^+(\HhL)$. Since for every $W \in \Repl(\h|_\p)$ the $k$-vector space $\Hhn{W}$ is a quotient of $\Ff^+(\HhL)$, also the latter is equipped with a filtration $F^*\Hhn{W}$ inducing the associated decomposition $\Hhn{W}=\bigoplus_{d \leq 0} \Hhn{W}(d)$ where $\Hhn{W}(d)=F^{d} \Hhn{W} \Qdx F^{d-1} \Hhn{W}$.

\begin{rmk} Once we choose local coordinates and an isomorphism between $\hL$ and $\gL$ we observe that the elements of $\Ff^+(\HhL)$ are $k[(c+\hslash)^{-1}]$-linear combinations of elements $X_r t^{-k_r} \circ \cdots \circ X_1t^{-k_1} \circ e_0$ with $k_r \geq \cdots \geq k_1 \geq 0$ and $r \geq 0$, where $e_0$ stands for $1 \in k$. We can explicitly write the graded pieces of $\Ff^+\HhL$ as
\[\Ff^+(\HhL)(-d)= \left\langle X_r t^{-k_r} \circ \cdots \circ X_1t^{-k_1} \circ e_0 \, | \, \sum_{i=1}^r k_i=d \right\rangle
\]
so that it is not zero only for $d\leq 0$ and in particular $\Ff^+(\HhL)(0)=k $ which shows that $\Hhn{W}(0)=W$.
\end{rmk}

The key ingredient to provide a morphism between $\VV_\ell(\V)_\Xx$ and $\VV_\ell(V)_X[\![\tau]\!]$ lies in the construction of the element $\epsilon(W)$ given by the following Proposition which we can see as a consequence of \cite[Lemma 6.5]{looijenga2005conformal}.

\begin{prop} \label{prop-epsilon} Let $W \in \Repl(\h|_\p)$ and $b_W^0 \colon W \otimes W^* \to k$ be the trace morphism. Then there exists an element
\[ \epsilon(W) = \sum_{d \geq 0} \epsilon(W)_d \cdot \tau^d \in \left(\Hhn{W} \otimes \Hhn{W^*}\right)[\![\tau]\!]
\]
satisfying the following conditions:
\begin{enumerate}[label=(\alph*)]
\item \label{epsdualb} the constant term $\epsilon(W)_0 \in W \otimes W^*$ is the dual of $b_W^0$ and for every $d \in \ZZ_{\geq 0}$ we have $\epsilon(W)_d \in \Hhn{W}(-d) \otimes \Hhn{W^*}(-d)$;
\item $\epsilon(W)$ is annihilated by the image of $\h_{k[\![\tp, \tm]\!]}$ in $\overline{U} \widehat{\h}_{\Ll_N}[\![\tau]\!]$.
\end{enumerate}
\end{prop}

\proof We choose an isomorphism between $\hLi[N]$ and $\g \otimes \Ll_N$, as well as an isomorphism between $\h_{\Op_\p}$ and $\g \otimes \Op_p$. The construction of $\epsilon(W)$ essentially lies in showing that the pairing $b_W^{(0)} \colon W \otimes W^* \to k$ extends to a unique pairing \[b_{W} \colon \Hhn{W} \otimes \Hhn{W^*} \to k\] such that for all  $(u,v) \in \Hhn{W} \otimes \Hhn{W^*}$ we have 
\begin{equation}\label{eqtlem}
b_W(X\tp^m u, v)+b_W(u,X\tm^{-m}v)=0
\end{equation}
for all $m \in \ZZ$ and $X \in \g$ and that $b_W$ is identically zero when restricted to $\Hhn{W}(d) \otimes \Hhn{W^*}(d')$ if $d \neq d'$. This is essentially \cite[Claim 1 of the proof of Proposition 6.2.1]{tsuchiya1989conformal} to which we refer. \endproof

We saw how to attach to any representation $W$ the element $\epsilon(W)$: we now use these elements to obtain the isomorphism map between $\VV_\ell(\V)_\Xx$ and $\VV_\ell(\V_0)_X[\![\tau]\!]$. The following statement, combined with Proposition \ref{prop-NodalDegenerationH} implies Theorem \ref{thm-HhnTrivDef}.

\begin{prop} \label{EisIso} The $k[\![\tau]\!]$ linear map 
\begin{align*} E \colon \Hhn{\V} \subset \Hhn{\V_0}[\![\tau]\!] &\to \bigoplus_{W \in \Repl(\h|_\p)} \Hhn{\V_0, W,  W^*}[\![\tau]\!]\\
u=\sum_{i\geq 0} u_i \tau^i &\mapsto \left( u \otimes \epsilon(W)\right)_{W \in \Repl(\h|_\p)} = \left(\sum_{i,d\geq 0} u_i \otimes \epsilon(W)_d\tau^{i+d}\right)_{W \in \Repl(\h|_\p)}
\end{align*}
induces the isomorphism
\[E_\hA \colon \VV_\ell(\V)_\Xx \to \bigoplus_{W \in \Repl(\h|_\p)} \Hhn{\V_0, W,  W^*}_{X_N}[\![\tau]\!].\]
of $k[\![\tau]\!]$-modules.
\end{prop}

\proof In order to prove that $E_{\hA}$ is an isomorphism we first quotient out by $\tau$ and using the identifications observed in Remark \ref{rmk-modulotau} we get the map
\[\![E_\hA]_{\tau=0} \colon \h_A \Qsx \Hhn{\V_0} \to \bigoplus_{W \in \Repl(\h|_\p)} \h_{{A_N}^*} \Qsx \Hhn{\V_0, W, W^*}\]
which sends the class of $u$ to $\left(u_0 \otimes \epsilon(W)_0 \right)_{W \in \Repl(\h|_\p)}$. Property \ref{epsdualb} of $\epsilon(W)_0$ tells us that $[E_\hA]_{\tau=0}$ is, up to some invertible factors, the inverse of the morphism induced by the $\{b_W\}$, which we showed to be an isomorphism in Proposition \ref{prop-NodalDegenerationH}.
Since $\hA \Qsx \Hhn{\V}$ is finitely generated, Nakayama's lemma and the fact that the right hand side is a free $k[\![\tau]\!]$-module guarantee that $E_\hA$ is an isomorphism.
\endproof

\begin{rmk} \label{rmk-locfreefamily} The argument we used run similarly if instead of starting with a covering of curves over $\Spec(k)$, we had considered a family of coverings $(\Xt \overset{q}\to X \overset{\pi}\to S, \sigma)$ where the singular locus of $X$ is given by one (or more) sections of $\pi$ and whose normalization is a covering of versal pointed smooth curves. Using these assumptions we are able to construct the canonical smoothing $\Xxt \to \Xx$ of $\Xt \to X$ over $S[\![\tau]\!]$ which is moreover a versal deformation of $(\Xt \overset{q}\to X \overset{\pi}\to S, \sigma)$. Once we have this construction, the analogue of Theorem \ref{thm-HhnTrivDef} follows. 
\end{rmk}

We have now all the ingredients to prove Theorem \ref{thm-LocFree} which, we recall, stated that $\VV_\ell(\V_1, \dots, \V_n)$ is locally free on the whole $\bHur[n]$.

\proof[Proof of Theorem \ref{thm-LocFree}] Let consider only the case $n=1$. From Corollary \ref{cor-ConnectionLocFree} we already know that $\VV_\ell(\V)$ is locally free on $\Hur$, so we have to check that this property extends to the boundary. Let $(\Xt \overset{q}\to X \to \Spec(k), \sigma_0)$ be a $k$-point of $\bHur \setminus \Hur$. As already mentioned we are left to show that $\VV_\ell(\V)$ is locally free on a neighbourhood of $(\Xt \overset{q}\to X \to \Spec(k), \sigma_0)$, i.e. that for one (hence any) versal deformation $(\Xxt \to \Xx \to S, \sigma)$ of $(\Xt \to X \to \Spec(k), \sigma_0)$, the $\Oo_S$-module $\VV_\ell(\V)_\Xx$ is locally free. Assume, for simplicity only, that $\p \in X(k)$ is the only nodal point of $X$. Consider the normalization $(\Xt_N \to X_N, \sigma_0, \p_+,\p_-)$ of $(\Xt \to X, \sigma_0)$ and denote by $(\Xxt_N \to \Xx_N \to S, \sigma, \mathfrak{P}_+, \mathfrak{P}_-)$ its universal deformation. Since we can see $\Xt \to X$ as a fibre of the covering obtained from $\Xxt_N \to \Xx_N$ by identifying $\mathfrak{P}_-$ and $\mathfrak{P}_+$, the previous remark allows us to conclude.\endproof

\bigskip
\small
\noindent\textit{Acknowledgments.} The main results of this paper are part of my PhD thesis, which was written in 2017 at the Universität Duisburg-Essen under the supervision of Jochen Heinloth. I am indebted to him for the constant support, for the useful discussions and comments on a preliminary version of this manuscript. Many thanks to Christian Pauly and Angela Gibney. I am grateful to the anonymous referee for their comments and suggestions.

\normalsize
\bibliography{Biblio}

\begin{thebibliography}{10}

\bibitem{ArbarelloCornalba}
{\sc E.~Arbarello, M.~Cornalba, and P.~A. Griffiths}, {\em Geometry of
  algebraic curves. {V}olume {II}}, vol.~268 of Grundlehren der Mathematischen
  Wissenschaften [Fundamental Principles of Mathematical Sciences], Springer,
  Heidelberg, 2011.
\newblock With a contribution by Joseph Daniel Harris.

\bibitem{balaji2011moduli}
{\sc V.~Balaji and C.~S. Seshadri}, {\em Moduli of parahoric
  {$\mathscr{G}$}-torsors on a compact {R}iemann surface}, J. Algebraic Geom.,
  24 (2015), pp.~1--49.

\bibitem{beauville1994conformal}
{\sc A.~Beauville}, {\em Conformal blocks, fusion rules and the {V}erlinde
  formula}, in Proceedings of the {H}irzebruch 65 {C}onference on {A}lgebraic
  {G}eometry ({R}amat {G}an, 1993), vol.~9 of Israel Math. Conf. Proc.,
  Bar-Ilan Univ., Ramat Gan, 1996, pp.~75--96.

\bibitem{BeauvilleLaszlo1994Conformal}
{\sc A.~Beauville and Y.~Laszlo}, {\em Conformal blocks and generalized theta
  functions}, Comm. Math. Phys., 164 (1994), pp.~385--419.

\bibitem{BGM15Vanishing}
{\sc P.~Belkale, A.~Gibney, and S.~Mukhopadhyay}, {\em Vanishing and identities
  of conformal blocks divisors}, Algebr. Geom., 2 (2015), pp.~62--90.

\bibitem{BGM16Nonvanishing}
{\sc P.~Belkale, A.~Gibney, and S.~Mukhopadhyay}, {\em Nonvanishing of
  conformal blocks divisors on {$\overline M_{0,n}$}}, Transform. Groups, 21
  (2016), pp.~329--353.

\bibitem{bertin2007champs}
{\sc J.~Bertin and M.~Romagny}, {\em Champs de {H}urwitz}, M\'em. Soc. Math.
  Fr. (N.S.),  (2011), p.~219.

\bibitem{bosch1990neron}
{\sc S.~Bosch, W.~L\"utkebohmert, and M.~Raynaud}, {\em N\'eron models},
  vol.~21 of Ergebnisse der Mathematik und ihrer Grenzgebiete (3) [Results in
  Mathematics and Related Areas (3)], Springer-Verlag, Berlin, 1990.

\bibitem{DrinfeldSimpson}
{\sc V.~G. Drinfel$\prime$d and C.~Simpson}, {\em {$B$}-structures on
  {$G$}-bundles and local triviality}, Math. Res. Lett., 2 (1995),
  pp.~823--829.

\bibitem{Edixhoven1992Neron}
{\sc B.~Edixhoven}, {\em N\'eron models and tame ramification}, Compositio
  Math., 81 (1992), pp.~291--306.

\bibitem{Fak12CBdivisors}
{\sc N.~Fakhruddin}, {\em Chern classes of conformal blocks}, in Compact moduli
  spaces and vector bundles, vol.~564 of Contemp. Math., Amer. Math. Soc.,
  Providence, RI, 2012, pp.~145--176.

\bibitem{Faltings1994Verlinde}
{\sc G.~Faltings}, {\em A proof for the {V}erlinde formula}, J. Algebraic
  Geom., 3 (1994), pp.~347--374.

\bibitem{EGAIII2}
{\sc A.~Grothendieck}, {\em {\'El\'ements de g\'eom\'etrie alg\'ebrique. {III}.
  \'Etude cohomologique des faisceaux coh\'erents. {II}}}, Inst. Hautes
  \'Etudes Sci. Publ. Math.,  (1963), p.~91.

\bibitem{heinloth2010uniformization}
{\sc J.~Heinloth}, {\em Uniformization of {$\mathscr{G}$}-bundles}, Math. Ann.,
  347 (2010), pp.~499--528.

\bibitem{KumarHong}
{\sc J.~Hong and S.~Kumar}, {\em Conformal blocks for galois covers of
  algebraic curves}, arXiv preprint arXiv:1807.00118,  (2018).

\bibitem{kac1994infinite}
{\sc V.~G. Kac}, {\em Infinite-dimensional {L}ie algebras}, Cambridge
  University Press, Cambridge, third~ed., 1990.

\bibitem{rainakac1988lecture}
{\sc V.~G. Kac and A.~K. Raina}, {\em Bombay lectures on highest weight
  representations of infinite-dimensional {L}ie algebras}, vol.~2 of Advanced
  Series in Mathematical Physics, World Scientific Publishing Co., Inc.,
  Teaneck, NJ, 1987.

\bibitem{Kumar}
{\sc S.~Kumar}, {\em Demazure character formula in arbitrary {K}ac-{M}oody
  setting}, Invent. Math., 89 (1987), pp.~395--423.

\bibitem{LS1997PicardBunG}
{\sc Y.~Laszlo and C.~Sorger}, {\em The line bundles on the moduli of parabolic
  {$G$}-bundles over curves and their sections}, Ann. Sci. \'Ecole Norm. Sup.
  (4), 30 (1997), pp.~499--525.

\bibitem{looijenga2005conformal}
{\sc E.~Looijenga}, {\em From {WZW} models to modular functors}, in Handbook of
  moduli. {V}ol. {II}, vol.~25 of Adv. Lect. Math. (ALM), Int. Press,
  Somerville, MA, 2013, pp.~427--466.

\bibitem{Mathieu}
{\sc O.~Mathieu}, {\em Formules de caract\`eres pour les alg\`ebres de
  {K}ac-{M}oody g\'en\'erales}, Ast\'erisque,  (1988), p.~267.

\bibitem{PappasRapoport2007Questions}
{\sc G.~Pappas and M.~Rapoport}, {\em Some questions about
  {$\mathscr{G}$}-bundles on curves}, in Algebraic and arithmetic structures of
  moduli spaces ({S}apporo 2007), vol.~58 of Adv. Stud. Pure Math., Math. Soc.
  Japan, Tokyo, 2010, pp.~159--171.

\bibitem{Pauly1996Parabolic}
{\sc C.~Pauly}, {\em Espaces de modules de fibr\'es paraboliques et blocs
  conformes}, Duke Math. J., 84 (1996), pp.~217--235.

\bibitem{sorger1996formule}
{\sc C.~Sorger}, {\em La formule de {V}erlinde}, Ast\'erisque,  (1996),
  pp.~Exp.\ No.\ 794, 3, 87--114.
\newblock S\'eminaire Bourbaki, Vol. 1994/95.

\bibitem{tsuchiya1989conformal}
{\sc A.~Tsuchiya, K.~Ueno, and Y.~Yamada}, {\em Conformal field theory on
  universal family of stable curves with gauge symmetries}, in Integrable
  systems in quantum field theory and statistical mechanics, vol.~19 of Adv.
  Stud. Pure Math., Academic Press, Boston, MA, 1989, pp.~459--566.

\bibitem{zelaci2017moduli}
{\sc H.~Zelaci}, {\em Moduli spaces of anti-invariant vector bundles and
  twisted conformal blocks}, arXiv preprint arXiv:1711.08296,  (2017).

\end{thebibliography}
\bibliographystyle{siam}
  
    \vfill 

\end{document}